\documentclass[fontsize=10pt]{scrreprt}

\usepackage[square,numbers,nonamebreak]{natbib}
\usepackage{setspace}

\usepackage[english]{babel}

\usepackage{chngcntr}

\usepackage[hang]{footmisc}
 at 11truept
\font\ssc=pplrc9d at 11 truept
\newcommand\qedbox{$\rlap{$\sqcap$}\sqcup$}

\usepackage[hmarginratio=1:1, bottom=2cm, top=2.5cm]{geometry}

\usepackage[automark]{scrlayer-scrpage}
\cohead{\textnormal{Joins of $\sigma$-subnormal subgroups}}
\lehead{\pagemark}
\rohead{\pagemark}

\let\ceheadL\cehead
\renewcommand\cehead[1]{
\ceheadL{\textnormal{#1}}
}

\cehead{Maria Ferrara --- Marco Trombetti}
\lefoot{}
\rofoot{}

\usepackage{graphicx}
\usepackage[usenames,dvipsnames,svgnames,table]{xcolor}
\definecolor{Maroon}{cmyk}{0, 0.87, 0.68, 0.32}
\definecolor{RoyalBlue2}{cmyk}{80,100,0,0.1}

\newcommand\auths[1]{\large \textsc{\textcolor{Maroon}{#1}}\setstretch{1.2}}
\newcommand\titl[1]{\center \linespread{1.1}\color{RoyalBlue2}\Large\textbf{ #1}\color{black}\bigskip} 
\renewcommand\abstract[1]{
\begin{center}
{\textbf{Abstract}}
\end{center}
{
\linespread{1.1}\fontsize{9pt}{-10pt}\selectfont #1}}

\usepackage[utf8]{inputenc}
\usepackage[sc,osf]{mathpazo}
\usepackage[euler-digits]{eulervm}
\usepackage[T1]{fontenc}
\usepackage{mathtools}

\DeclareSymbolFont{operators}{\encodingdefault}{ppl}{m}{n}
\DeclareMathAlphabet{\mathbf}{\encodingdefault}{ppl}{bx}{n}
\DeclareMathAlphabet{\mathit}{\encodingdefault}{ppl}{m}{it}

\usepackage{amsfonts}
\usepackage{amsmath}
\usepackage{ragged2e}
\usepackage{titlesec}
\usepackage{amssymb}

\renewcommand{\thesection}{\arabic{section}}
 
\titleformat{\section}{\medskip\bigskip\normalfont\Large\bf}{\thesection}{0.5em}{}
\titleformat{\subsection}{\smallskip\bigskip\normalfont\large\bf}{\thesubsection}{0.5em}{}

\usepackage{amsthm}
\newtheoremstyle{dotless}{}{}{\itshape}{}{\bfseries}{}{1em}{}

\theoremstyle{dotless}
\newtheorem{theo}{Theorem}

\newtheorem{lem}[theo]{Lemma}
\newtheorem{cor}[theo]{Corollary}

\newtheorem{quest}[theo]{Question}
\renewenvironment{proof}{\smallbreak\noindent {\sc Proof \;---\;}}{\hfill\qedbox}

\counterwithout{theo}{section}
\numberwithin{theo}{section}

\DeclareOldFontCommand{\rm}{\normalfont\rmfamily}{\mathrm}
\DeclareOldFontCommand{\sf}{\normalfont\sffamily}{\mathsf}
\DeclareOldFontCommand{\tt}{\normalfont\ttfamily}{\mathtt}
\DeclareOldFontCommand{\bf}{\normalfont\bfseries}{\mathbf}
\DeclareOldFontCommand{\it}{\normalfont\itshape}{\mathit}
\DeclareOldFontCommand{\sl}{\normalfont\slshape}{\@nomath\sl}
\DeclareOldFontCommand{\sc}{\normalfont\scshape}{\@nomath\sc}

\usepackage{tikz}

\begin{document}

\setheadsepline{1pt}[\color{black}]

\titl{Joins of $\sigma$-Subnormal Subgroups
\footnote{The authors are supported by GNSAGA (INdAM) and are members of the non-profit association ``AGTA --- Advances in Group Theory and Applications'' (www.advgrouptheory.com). }}

\auths{Maria Ferrara --- Marco Trombetti}

\thispagestyle{empty}
\justify\noindent
\setstretch{0.3}
\abstract{Let $\sigma=\{\sigma_j\,:\, j\in J\}$ be a partition of the set $\mathbb{P}$ of all prime numbers. A subgroup $X$ of a finite group $G$ is~{\it $\sigma$-subnormal} in $G$ if there exists a chain of subgroups $$X=X_0\leq X_1\leq\ldots\leq X_n=G$$ such that, for each $1\leq i\leq n-1$, $X_{i-1}\trianglelefteq X_i$ or $X_i/(X_{i-1})_{X_i}$ is a $\sigma_{j_i}$-group for some $j_i\in J$. Skiba~\cite{skiba} studied the main properties of~\hbox{$\sigma$-sub}\-normal subgroups in finite groups and showed that the set of all~\hbox{$\sigma$-sub}\-normal subgroups plays a very relevant role in the structure of a finite soluble group. In \cite{MFsigma}, we laid the foundation of a general theory of $\sigma$-subnormal subgroups (and $\sigma$-series) in locally finite groups. It turns out that the main difference between the finite and the locally finite case concerns the behaviour of the join of $\sigma$-subnormal subgroups: in finite groups, $\sigma$-subnormal subgroups form a sublattice of the lattice of all subgroups (see for instance \cite{Aguilera}), but this is no longer true for arbitrary locally finite groups. This situation is very similar to that concerning with subnormal subgroups, so, as in the case of subnormal subgroups (see for instance \cite{Rob65},\cite{mathz},\cite{Smith},\cite{williams} and the monograph \cite{lennox}), it therefore makes sense to study the class $\mathfrak{S}_\sigma^\infty$ (resp. $\mathfrak{S}_\sigma$) of locally finite groups in which the join of (resp. finitely many)~\hbox{$\sigma$-sub}\-normal subgroups is $\sigma$-subnormal. In particular, the aim of this paper is to study how much one can extend a group in one of these classes before going outside the same class (see for example Theorems~\ref{sinftymaxsn}, \ref{sinftyminmax}, \ref{theosopra} and \ref{theorankunico}). Furthermore, some $\sigma$-subnormality criteria for the join of two $\sigma$-subnormal subgroups are obtained: for example, similarly to a celebrated theorem of Williams (see \cite{williams}), we give a necessary and sufficient condition for a join of two $\sigma$-subnormal subgroups to always be $\sigma$-subnormal; as a consequence we show that the join of two orthogonal $\sigma$-subnormal subgroups is $\sigma$-subnormal (this is the analog of a result of Roseblade \cite{mathz}).}

\setstretch{2.1}
\noindent
{\fontsize{10pt}{-10pt}\selectfont {\it Mathematics Subject Classification \textnormal(2020\textnormal)}: 20F50, 20E15}\\[-0.8cm]

\noindent 
\fontsize{10pt}{-10pt}\selectfont  {\it Keywords}: locally finite group; $\sigma$-subnormal subgroup; subnormal subgroup\\[-0.8cm]

\setstretch{1.1}
\fontsize{11pt}{12pt}\selectfont
\section{Introduction}

Subnormal subgroups were introduced by Helmut Wielandt in an attempt to bypass the fact that normality is not transitive in general, and nowadays they are an invaluable tool in studying both finite and infinite groups; the reader is referred to~\cite{lennox}, where the main properties and some of the applications of subnormal subgroups are exhibited. A well known theorem of~Wielandt shows that any join of subnormal subgroups in a finite group is subnormal, but this is no longer the case for infinite groups (there exist even locally finite, metabelian counterexamples). A large part of the literature in this area is in fact devoted to the study of criteria for which the join of (two) subnormal subgroups is subnormal (see for instance \cite{mathz2}, \cite{mathz}, \cite{stone}, \ldots) and to the study of the classes of groups in which the joins of (two) subnormal subgroups are always subnormal (see for instance \cite{Rob65}, \cite{Smith}, \cite{williams}, \ldots). There are many beautiful results here, and we now quote some of them that are relevant to our paper. Let $H$ and~$K$ be subnormal subgroups of a group $G$. It turns out that if $HK=KH$, then $\langle H,K\rangle$ is subnormal in~$G$ (see \cite{lennox}, Theorem 1.2.5). Also, a  theorem of Roseblade asserts that if $H$ and~$K$ are {\it orthogonal} (\hbox{i.e.} $H/H'\otimes K/K'=\{0\}$), then $\langle H,K\rangle$ is subnormal in $G$ (see~\cite{mathz}). This is in fact a very special case of a celebrated theorem of Williams \cite{williams} giving necessary and sufficient conditions for the join of two subnormal subgroups to always be subnormal. Furthermore, we recall that $J=\langle H,K\rangle$ is subnormal in~$G$ provided that $J'$ satisfies the maximal condition on subgroups (see~\cite{lennox}, Corollary~1.3.9) or has finite rank (see~\cite{lennox},~The\-o\-rem~3.4.3). Now, let $\mathfrak{S}^\infty$ (resp. $\mathfrak{S}$) be the class of groups in which the join of (resp. two) subnormal subgroups is subnormal. It has been proved in \cite{Smith} that if $G$ is any group having normal subgroups $N\leq M$ such that $M'N/N$ and~$G/M$ have {\it finite rank} (a group is said to have {\it finite rank} $r$ if any of its finitely generated subgroups can be generated by $r$ elements), and $N\in\mathfrak{S}^\infty$, then $G$ belongs to $\mathfrak{S}$; moreover, any extension of a group of finite rank by an $\mathfrak{S}$-group is an $\mathfrak{S}$-group. Also, if~$\mathfrak{M}$ denotes the class of groups having a finite series whose factors satisfy the maximal or the minimal condition on subnormal subgroups, then $\mathfrak{MS}^\infty=\mathfrak{S}^\infty$ (see~\cite{lennox}, Theorem~2.6.1); here, if $\mathfrak{X}$ and $\mathfrak{Y}$ are classes of groups, we denote by $\mathfrak{XY}$ the class of all~\hbox{\it $\mathfrak{X}$-by-$\mathfrak{Y}$ groups}, \hbox{i.e.} all groups $G$ having a normal subgroup $N\in\mathfrak{X}$ such that~\hbox{$G/N\in\mathfrak{Y}$.}

\medskip

Recently, Skiba \cite{skiba} began the study of a handful and powerful generalization of subnormal subgroups in finite groups. If $\sigma$ is a partition of the set $\mathbb{P}$ of all prime numbers, $G$ is a finite group and $X$ is a subgroup of $G$, then he said that $X$ is~{\it $\sigma$-sub\-normal} in $G$ if there is a chain of subgroups $$X=X_0\leq X_1\leq\ldots\leq X_n=G$$ such that, for each $1\leq i\leq n-1$, $X_{i-1}\trianglelefteq X_i$ or $X_i/(X_{i-1})_{X_i}$ is a $\sigma_{j_i}$-group for some~\hbox{$\sigma_{j_i}\in\sigma$.} Skiba showed that the set of all $\sigma$-subnormal subgroups has a strong influence on the structure of a finite soluble group and this led many authors to investigate which of the most relevant theorems about subnormal subgroups have analogs in terms of $\sigma$-subnormal subgroups (see for instance \cite{adolfo}, \cite{perez}, \cite{Aguilera}). It turns out that one of the main features of $\sigma$-subnormal subgroups (in finite groups) is that the join of two $\sigma$-subnormal subgroups is always $\sigma$-subnormal (see for instance \cite{Aguilera}), so they form a sublattice of the lattice of all subgroups. 

\medskip

In \cite{MFsigma}, we laid a foundation of a general theory of $\sigma$-subnormal subgroups (and \hbox{$\sigma$-series)} in {\it locally finite groups} (\hbox{i.e.} groups in which every finitely generated subgroup is finite) and we showed that, although the join of two $\sigma$-subnormal subgroups is not always $\sigma$-subnormal, this kind of subgroups (and its generalizations) often behaves better than subnormal subgroups. As far as the join is concerned, it has been proved for instance that the join of two $\sigma$-subnormal subgroups $H$ and $K$ of a locally finite group is $\sigma$-subnormal provided that $HK=KH$ (see \cite{MFsigma}, Theorem 3.15); and that the join of two $\sigma$-subnormal subgroups of a locally finite group $G$ is always $\sigma$-subnormal when $G$ is linear (see \cite{MFsigma}, Theorem 3.35) or when $G$ satisfies the maximal condition on subnormal subgroups (see \cite{MFsigma}, Theorem 3.25).

\medskip

The aim of this paper is to study the class $\mathfrak{S}_\sigma^\infty$ (resp. $\mathfrak{S}_\sigma$) of all locally finite groups in which the join of (resp. finitely many) $\sigma$-subnormal subgroups is $\sigma$-subnormal. We show that, although these classes are not extension closed, we can extend them through the class $\mathfrak{M}$. In this respect, the main results can be summarized as follows:

\begin{itemize}

\item[(i)] $\mathfrak{MS}_\sigma^\infty\mathfrak{M}=\mathfrak{S}_\sigma^\infty$ (see Theorems \ref{sinftymaxsn} and \ref{sinftyminmax});

\item[(ii)] $\mathfrak{MS}_\sigma\mathfrak{M}=\mathfrak{S}_\sigma$ (see Theorems \ref{sinftyminmaxfacile} and \ref{theosopra});

\item[(iii)] $\mathfrak{S}_\sigma^\infty\mathfrak{Y}\leq\mathfrak{S}_\sigma$ (see Theorem \ref{dacitare}), where $\mathfrak{Y}$ is the class of all locally finite groups having a normal subgroup $N$ such that $G/N$ and $N'$ have finite rank;

\item[(iv)] $\mathfrak{S}_\sigma^\infty\mathfrak{T}\leq\mathfrak{S}_\sigma$ (see Theorem \ref{corabelianinftheo}), where $\mathfrak{T}$ is the class of all locally finite groups in which normality is a transitive relation.

\end{itemize}

\medskip

It is also interesting to observe that if the partition $\sigma$ is finite, then the join of two~\hbox{$\sigma$-sub}\-normal subgroups is always $\sigma$-subnormal (see Theorem \ref{finalone}) although the corresponding statement for the join of arbitrarily many~\hbox{$\sigma$-sub}\-normal subgroups is still in doubt (see Question \ref{questionrelevant}).

As a by-product of the study of the class $\mathfrak{S}_\sigma$, we are able to provide many~\hbox{$\sigma$-sub}\-normality criteria for the join of two~\hbox{$\sigma$-sub}\-normal subgroups. For example, we are able to show that the theorems of Williams and Roseblade extend to~\hbox{$\sigma$-sub}\-normal subgroups (see Theorem \ref{williamstheo} and Corollary \ref{roseblade}) and that the join of two~\hbox{$\sigma$-sub}\-normal subgroups is $\sigma$-subnormal provided that the commutator subgroup is in $\mathfrak{M}$ (see Theorem \ref{corderivato}) or has finite rank (see the remark before Theorem \ref{finalone}). These results have been obtained through the study of the subnormality properties of the $\sigma$-hypercentral residuals (see Theorem \ref{theopermuta}); in this respect, we also quote two further interesting consequences, namely,~The\-o\-rem~\ref{permutizer} and~The\-o\-rem~\ref{altroteorem}.

\medskip

Basic notation, definitions and preliminary results are dealt with in Section \ref{preliminary}. Unexplained notation is intended to be standard and their meaning can be found in~\cite{Rob72} or~\cite{lennox}. In Section \ref{sec3} and Section \ref{sec4} we deal respectively with the class~$\mathfrak{S}_\sigma^\infty$ and the class~$\mathfrak{S}_\sigma$.  Actually, in Section \ref{sec3} we also need to deal with the class $\mathfrak{S}_{\operatorname{sn}\rightarrow\sigma}^\infty$ of all locally finite groups in which the join of arbitrarily many subnormal subgroups is~\hbox{$\sigma$-sub}\-normal; note that by Theorem \ref{corquestion2} , the class $\mathfrak{S}_\sigma$ coincides with the class of locally finite groups in which the join of two subnormal subgroups is $\sigma$-subnormal. In Section \ref{lastsection}, the main $\sigma$-subnormality criteria are deduced.

A diagram showing the relationship among all the classes of groups we study in this paper can be found at the beginning of Section \ref{sec4}.

\section{Preliminaries and notation}\label{preliminary}

In this section we deal with the basic definitions and results concerning $\sigma$-subnormal subgroups and related concepts. In the whole paper, if not otherwise stated, the Greek letter $\sigma$ is always meant to denote a partition of the set $\mathbb{P}$ of all prime numbers.

\subsection*{$\sigma$-Nilpotency and related concepts}

Let $\tau$ be any set of (disjoint) subsets of $\mathbb{P}$. A group $G$ is~\hbox{{\it $\tau$-primary}} (or a {\it $\tau$-group}) if there is~\hbox{$\tau_j\in\tau$} such that $\pi(G)\subseteq\tau_j$; here, if $S$ is a subset of $G$, the symbol $\pi(S)$ denotes the set of all primes dividing at least one of the orders of the elements of $S$. Moreover, a group is {\it $\tau$-soluble} (resp. {\it hyper-$\tau$}) if it has a finite (resp. ascendant) series whose factors are~\hbox{$\tau$-pri}\-mary; the smallest length of such a series is the~{\it $\tau$-length} of $G$ (in the case of hyper-$\tau$ groups, the length is of course an ordinal number). Remark that a~\hbox{$\tau$-soluble} group is hyper-$\tau$, and that a hyper-$\tau$ group is periodic;  also note that in the definition of ‘‘hyper-$\tau$’’ we can ask for the series to be normal and even characteristic (replace the first term by the product of all normal~\hbox{$\tau_j$-sub}\-groups for a suitable $\tau_j\in\tau$, and then continue by transfinite recurrence). Now, a group $G$ is {\it $\tau$-hypercentral} if it is a direct product of $\tau$-groups; in this case, $G$ is also called {\it $\tau$-nilpotent} if $\pi(G)$ is covered by finitely many elements of $\tau$. 

Let $G$ be a group and $\tau_j\in\tau$. Then the unique maximal normal $\tau_j$-subgroup of $G$ is the {\it $\tau_j$-com\-po\-nent} of $G$. A {\it $\tau$-component} of $G$ is just one of its $\tau_i$-components for some $\tau_i\in\tau$. Thus, for example, if $\sigma$ is a partition of $\mathbb{P}$, then a group is $\sigma$-hypercentral if and only if it is the product of its $\sigma$-components.

Finally, we remark that, for any partition $\sigma$ of $\mathbb{P}$, the property of being $\sigma$-hyper\-central can be detected by the behaviour of the finitely generated subgroups, so is a local property.

\subsection*{$\sigma$-Subnormality}

Let $\sigma=\{\sigma_i\,:\, i\in I\}$ be an arbitrary partition of the set $\mathbb{P}$ of all prime numbers. Let~$(\mathcal{I},<)$ be a linearly ordered set. Let $G$ be a group and $X\leq G$. A {\it $\sigma$-series} between~$X$ and~$G$ is a set $\mathcal{S}=\big\{(\Lambda_j,V_j)\,:\, j\in \mathcal{I}\big\}$ such that:

\begin{itemize}
\item[(i)] for each $j\in \mathcal{I}$, $X\leq V_j\leq \Lambda_j\leq G$;
\item[(ii)] $G\setminus X=\bigcup_{j\in \mathcal{I}}(\Lambda_j\setminus V_j)$;
\item[(iii)] $\Lambda_k\leq V_j$ if $k<j$;
\item[(iv)] for each $j\in\mathcal{I}$, either $V_j$ is normal in $\Lambda_j$, or $\Lambda_j/(V_j)_{\Lambda_j}$ is $\sigma$-primary.
\end{itemize}

\noindent In other words, a $\sigma$-series between $X$ and $G$ is a chain $\mathcal{S}$ of subgroups of $G$ such that either $K$ is normal in $H$ or $H/K_H$ is $\sigma$-primary, whenever $H<K$ are elements of $\mathcal{S}$ for which there is no element $L$ of $\mathcal{S}$ between them.

If there is a $\sigma$-series between $X$ and $G$, the subgroup $X$ is called {\it $\sigma$-serial} in $G$.  If the index set $\mathcal{I}$ is finite, then $X$ is {\it $\sigma$-subnormal} in $G$; while $X$ is {\it $\sigma$-ascendant} (resp. {\it $\sigma$-de\-scend\-ant}) if~$\mathcal{I}$ is well-ordered (resp. inversely well-ordered). Furthermore, for the sake of simplicity, we also say that $X$ is {\it $\sigma$-normal} in $G$, and we write $X\trianglelefteq_\sigma G$, whenever $|\mathcal{I}|=1$, or, in other words, if~$X$ is normal in $G$ or~$G/X_G$ is $\sigma$-primary. Thus, using this notation, a~\hbox{$\sigma$-sub}\-normal subgroup is just a subgroup $X$ for which there is a chain of subgroups of the following type \[\begin{array}{c}\label{star}
X=X_0\trianglelefteq_\sigma X_1\trianglelefteq_\sigma\ldots\trianglelefteq_\sigma X_n=G,\tag{$\star$}
\end{array}
\] and, in analogy with subnormality (recall that $X\trianglelefteq\trianglelefteq G$ means that $X$ is a subnormal subgroup of $G$), we may write $X\trianglelefteq\trianglelefteq_\sigma G$. If $\sigma_i\in\sigma$, we say that $X$ is {\it $\sigma_i$-normal} in $G$ if~$G/X_G$ is a $\sigma_i$-group; we denote this fact by writing $X\trianglelefteq_{\sigma_i} G$. Any $\sigma_i\in\sigma$ such that the expression $X_j\trianglelefteq_{\sigma_i} X_{j+1}$ appears in \eqref{star} is said to be {\it involved} in \eqref{star}.

It has been proved in \cite{MFsigma}, Theorem 3.1, that in locally finite groups, the homomorphic images of $\sigma$-serial subgroups are $\sigma$-serial (a fact that will be often employed in the paper without further notice), and that $\sigma$-serial subgroups play a very relevant role in understanding if a join of $\sigma$-subnormal subgroups is $\sigma$-subnormal. For example, it turns out that the join of $\sigma$-subnormal subgroups of a locally finite group is always $\sigma$-serial (see \cite{MFsigma}, Theorem 3.10) and that a $\sigma$-serial subgroup containing a $\sigma$-embedded subgroup which is subnormal in the whole group is~\hbox{$\sigma$-sub}\-normal~(see~\cite{MFsigma},~The\-o\-rem~3.6). Recall that a (normal) subgroup $H$ of a group~$G$ is~\hbox{$\sigma$\it -\textnormal{(}normally\textnormal{)} em}\-{\it bedded} in~$G$ if there is a finite (normal) series $$H=H_0\trianglelefteq H_1\trianglelefteq\ldots\trianglelefteq H_n=G$$ such that, for every $j=0,\ldots,n-1$, $H_{j+1}/H_j$ is a $\sigma_{i_j}$-group for some $\sigma_{i_j}\in\sigma$; in this case, we also say that $H$ is {\it $(\sigma_{i_0},\ldots,\sigma_{i_{n-1}})$-embedded in $G$}. In Lemma 2.9 of \cite{MFsigma}, we proved that every $\sigma$-subnormal subgroup of a group contains a characteristic~\hbox{$\sigma$-em}\-bedded subgroup that is subnormal in the whole group.

Finally, if $H$ is a subgroup of a locally finite group $G$ for which there is a $\sigma$-series $$H=H_0\trianglelefteq_{\sigma_1} H_1\trianglelefteq_{\sigma_2}\ldots\trianglelefteq_{\sigma_{n-1}} H_n=G,$$ we say that $H$ is {\it strictly $(\sigma_1,\ldots,\sigma_{n-1})$-subnormal} in $G$. It is easy to see that every such subgroup contains a characteristic $\sigma$-embedded subgroup (see \cite{MFsigma}, Lemma 2.10).

\medskip

The following lemma is a generalization (and a corollary) of Lemma 3.14 of~\cite{MFsigma}: it shows that there is no need to assume that both subgroups in that statement are~\hbox{$\sigma$-normal} in their join.

\begin{lem}\label{cor3.14}
Let $\sigma=\{\sigma_i\,:\, i\in I\}$ be a partition of $\mathbb{P}$. Let $G$ be a locally finite group. If~$H$ and~$K$ are $\sigma$-subnormal subgroups of $G$ and $H$ is $\sigma$-normal in $J=\langle H,K\rangle$, then $J$ is~\hbox{$\sigma$-sub}\-normal in $G$.
\end{lem}
\begin{proof}
By Theorem 3.15 of \cite{MFsigma} we may assume $H\trianglelefteq_{\sigma_i} J$ for some $i\in I$. Let $L$ be the normal core $H_J$ of $H$ in $J$; by definition, $J/L$ is a $\sigma_i$-group. Moreover, since $L$ is normal in $H$, it is $\sigma$-subnormal in $G$; hence  $KL$ is $\sigma$-subnormal in $G$ by Theorem 3.15 of \cite{MFsigma}. On the other hand, $(KL)_J$ contains $L$, so $KL$ is $\sigma$-normal in $G$. An application of~Lem\-ma~3.14 of \cite{MFsigma} yields that $J=\langle H, KL\rangle$ is $\sigma$-subnormal in $G$.
\end{proof}

\medskip

Finally, in connection with the previous concepts and definitions, we note the following two easy facts: (i) a locally finite group  is $\sigma$-hypercentral if and only if all its subgroups are $\sigma$-serial (use~Lemma~2.6 of \cite{MFsigma}); (ii) in a locally finite group $G$, all subgroups are $\sigma$-subnormal if and only if $G=H\times K$, where $\pi(H)\cap\pi(K)=\emptyset$, $H$ is~\hbox{$\sigma$-nil}\-potent, $K$ is $\sigma$-hypercentral, and all subgroups of $K$ are subnormal (use Lemmas 2.6 and 2.9 of \cite{MFsigma}).

\subsection*{Residuals}

Let $\mathfrak{X}$ be any class of groups. If $G$ is a group, we denote by $G^{\mathfrak{X}}$ the {\it $\mathfrak{X}$-residual} of~$G$, \hbox{i.e.} the intersection of all normal subgroups $N$ of $G$ such that $G/N$ is in $\mathfrak{X}$. If~$G$ is a locally finite group, $\pi\subseteq\mathbb{P}$, and $\mathfrak{X}$ is the class of all $\pi$-groups, then~$G^{\mathfrak{X}}$ is the~\hbox{\it $\pi$-residual} of~$G$, and is usually denoted by $O^{\pi}(G)$: it is the smallest normal subgroup $N$ of~$G$ such that~$G/N$ is a $\pi$-group. Moreover, if $\tau$ is any set of subsets of $\mathbb{P}$, we define the~\hbox{\it $\tau$-residual}~$G^\tau$ of a locally finite group~$G$ as the intersection of all~\hbox{$\pi$-residuals} of~$G$ for $\pi\in\tau$. Of course, if $\sigma$ is any partition of $\mathbb{P}$, then the $\sigma$-residual of a locally finite group is just its $\mathfrak{X}$-residual, where $\mathfrak{X}$ is the class of all $\sigma$-hypercentral groups. The~\hbox{$\sigma$-resi}\-duals of $\sigma$-subnormal subgroups play a very important part in the problem of the join of~\hbox{$\sigma$-sub}\-normal subgroups. First, it turns out that these are always subnormal in the whole group (this is a generalization of Lemma 2.9 of \cite{MFsigma}). It will be clear from the proof that the characteristic subgroup $N$ we find is such that $H/N$ is~$\{\sigma_{1},\ldots,\sigma_{m}\}$-nilpotent, where $\sigma_1,\ldots,\sigma_m$ are the elements of $\sigma$ that are involved in a given finite~\hbox{$\sigma$-series} connecting $H$ to $G$.

\begin{lem}\label{alternative}
Let $\sigma=\{\sigma_i\,:\, i\in I\}$ be a partition of $\mathbb{P}$. Let $H$ be a $\sigma$-subnormal subgroup of a locally finite group $G$. Then $H$ contains a characteristic subgroup $N$ such that $H/N$ is~\hbox{$\sigma$-nil}\-potent and $N\trianglelefteq\trianglelefteq G$.
\end{lem}
\begin{proof}
Let $$H=H_0\trianglelefteq_\sigma H_1\trianglelefteq_\sigma\ldots\trianglelefteq_\sigma H_n=G$$ be a $\sigma$-subnormal series connecting $H$ to $G$. The statement clearly holds if $n=0$. Assume $n>0$. By induction $H$ has a characteristic subgroup $M$ such that $H/M$ is~\hbox{$\sigma$-nil}\-potent and $M$ is subnormal in $H_{n-1}$. If $H_{n-1}\trianglelefteq H_n=G$, we are done. Assume $H_{n-1}\trianglelefteq_{\sigma_i} G$ for some $i\in I$. Then $G/K$ is a $\sigma_i$-group, where $K=(H_{n-1})_G$. Now, $H/M\cap K$ is $\sigma$-nilpotent and $M\cap K\trianglelefteq\trianglelefteq K\trianglelefteq G$, so the only thing left to do is to replace $M\cap K$ by the $\tau$-residual of $H$, where $\tau$ is the set of all $\sigma_j\in\sigma$ such that $\sigma_j\cap\pi(H/M\cap K)\neq\emptyset$.
\end{proof}

\medskip

Second, $\sigma$-residuals have very nice permutability properties. Using Lemma 3 (5) of~\cite{perez}, the argument in Lemma 3.3 of~\cite{Kamarnikov} proves the following result.

\begin{lem}\label{lemmakarmo}
Let $\sigma=\{\sigma_i\,:\, i\in I\}$ be a partition of $\mathbb{P}$ and let $G$ be a $\sigma$-soluble finite group. If~$H$ and $K$ are $\sigma$-subnormal subgroups of $G$, then $H^\sigma K=KH^\sigma$.
\end{lem}

\begin{cor}\label{karmo}
Let $\sigma=\{\sigma_i\,:\, i\in I\}$ be a partition of $\mathbb{P}$ and let $G$ be a $\sigma$-soluble finite group. If $H$ and $K$ are $\sigma$-subnormal subgroups of $G$, then $H^\sigma K^\sigma=K^\sigma H^\sigma$.
\end{cor}
\begin{proof}
Since $H^\sigma\trianglelefteq\trianglelefteq_\sigma G$, an application of Lemma \ref{lemmakarmo} to~\hbox{$\langle H^\sigma, K\rangle$} yields the result.~\end{proof}

\begin{theo}\label{theopermuta}
Let $\sigma=\{\sigma_i\,:\, i\in I\}$ be a partition of $\mathbb{P}$. Let $H$ and $K$ be $\sigma$-subnormal subgroups of a locally finite group $G$. Then $\langle H,K\rangle^\sigma=H^\sigma\, K^\sigma$. In particular, $H^\sigma K=KH^\sigma$ and the following conditions are equivalent:

\begin{itemize}
\item[\textnormal{(1)}] $HK=KH$;

\item[\textnormal{(2)}] $\langle H,K\rangle=HK\langle H,K\rangle^\sigma$;
\end{itemize}
\end{theo}
\begin{proof}
It is certainly possible to assume $G=\langle H,K\rangle$. Let $R,S$ and $T$ be the $\sigma$-soluble residuals of $H,K$ and $G$, respectively. 

Let $\mathcal{F}$ be the set of all finite subgroups of $G$ such that $F=\langle F\cap H,\, F\cap K\rangle$; clearly,~$\mathcal{F}$ is a local system of $G$. For each $F\in\mathcal{F}$, let $R(F),S(F)$ and $T(F)$ be the $\sigma$-soluble residuals of $F\cap H$, $F\cap K$ and $F$, respectively;  it follows from Corollary 6.5.48 of~\cite{ezquerro} that $T(F)=\langle R(F),\, S(F)\rangle$. Actually, $R(F)$ and $S(F)$ are perfect subnormal subgroups of $F$ by~Lem\-ma~\ref{alternative}, so $T(F)=R(F)\,S(F)$ by a well known theorem of Wielandt.

Let $$R_1=\bigcup_{F\in\mathcal{F}}R(F),\quad S_1=\bigcup_{F\in\mathcal{F}}S(F)\quad\textnormal{and}\quad T_1=\bigcup_{F\in\mathcal{F}}T(F).$$ Since $R(E)\geq R(F)$ where $E\geq F$ are in $\mathcal{F}$, it follows that $R_1$ is a normal subgroup of~$H$; similarly, $S_1$ and $T_1$ are normal subgroups of $K$ and $G$, respectively. Let $x\in R_1$ and $y\in S_1$. By definition, there are finite subgroups $E_1$ and $E_2$ in $\mathcal{F}$ such that $x\in R(E_1)$ and $y\in S(E_2)$. Choose $F\in\mathcal{F}$ such that $F\geq\langle E_1,E_2\rangle$. Since $R(F)\geq R(E_1)$, $S(F)\geq S(E_2)$, and $R(F)\, S(F)=S(F)\, R(F)$, we get that $xy=y'x'$, for some $y'\in R(F)\leq R_1$ and $x'\in S(F)\leq S_1$. It follows that $T_1=\langle R_1,S_1\rangle=R_1S_1$.

Clearly, $R_1\leq R\leq H^\sigma$ and $S_1\leq S\leq K^\sigma$, so if we can prove that $H^\sigma K^\sigma T_1=K^\sigma H^\sigma T_1$, then $$H^\sigma K^\sigma=H^\sigma R_1S_1K^\sigma=H^\sigma K^\sigma T_1=K^\sigma H^\sigma T_1=K^\sigma S_1R_1H^\sigma=K^\sigma H^\sigma.$$ Moreover, $T_1\leq T\leq G^\sigma$, so if we can prove that $(G/T_1)^\sigma=H^\sigma K^\sigma/T_1$, then $G^\sigma=H^\sigma K^\sigma$. Since $(HT_1/T_1)^\sigma=H^\sigma T_1/T_1$ and $(KT_1/T_1)^\sigma=K^\sigma T_1/T_1$, it is possible to assume~\hbox{$T_1=\{1\}$}. In particular,~$G$ is locally $\sigma$-soluble.

\smallskip

Now, in order to see that $G^\sigma=\langle H^\sigma,\, K^\sigma\rangle$, we only need to apply the above argument replacing $\sigma$-soluble residuals by $\sigma$-residuals, Corollary 6.5.48 of~\cite{ezquerro} by Lem\-ma~3~(5) of~\cite{perez}, and to notice that locally $\sigma$-hypercentral groups are $\sigma$-hypercentral (so in this case $R_1=H^\sigma$, and so on). The above argument also yields  that $H^\sigma K^\sigma =K^\sigma H^\sigma$: although the $\sigma$-residuals are not in general perfect, we may replace the theorem of~Wielandt by Corollary \ref{karmo}. 

Finally, $H^\sigma K=H^\sigma K^\sigma K=G^\sigma K=KG^\sigma=KK^\sigma H^\sigma =KH^\sigma$ and, if $G=HKG^\sigma$, then $G=HKG^\sigma=HH^\sigma K^\sigma K=HK$. The statement is proved.
\end{proof}

\medskip

The above result has many consequences. We now derive some of the most immediate among them and we leave the most relevant ones (dealing with the $\sigma$-subnormality of the join of two $\sigma$-subnormal subgroups) to Section \ref{lastsection}.

\begin{cor}
Let $\sigma=\{\sigma_i\,:\, i\in I\}$ be a partition of $\mathbb{P}$. Let $H$ and $K$ be $\sigma$-subnormal subgroups of a locally finite group $G$. If $H=H^\sigma$, then $HK=KH$.
\end{cor}

\begin{cor}\label{corollaryperm2}
Let $\sigma=\{\sigma_i\,:\, i\in I\}$ be a partition of $\mathbb{P}$. Let $H$ and $K$ be $\sigma$-subnormal subgroups of a locally finite group $G$. Then $\langle H,K\rangle^\sigma$ is subnormal in $G$.
\end{cor}
\begin{proof}
It follows from Theorem \ref{theopermuta} that $\langle H,K\rangle^\sigma=H^\sigma K^\sigma$. On the other hand,~Lem\-ma~\ref{alternative} shows that $H^\sigma$ and $K^\sigma$ are subnormal subgroups of $G$, and hence $\langle H,K\rangle^\sigma$ is subnormal in $G$, being the product of two subnormal subgroups.
\end{proof}

\begin{cor}\label{corfinitelymanysub}
Let $\sigma=\{\sigma_i\,:\, i\in I\}$ be a partition of $\mathbb{P}$. Let $\mathcal{H}=\{H_1,\ldots,H_n\}$ be a finite family of $\sigma$-subnormal subgroups of a locally finite group $G$. If $J=\langle H_1,\ldots, H_n\rangle$, then $J^\sigma=H_1^\sigma\ldots H_n^\sigma$.
\end{cor}
\begin{proof}
It follows from Theorem \ref{theopermuta} that $X=H_1^\sigma\ldots H_n^\sigma$ is a subgroup of $J$. In order to show that $X$ is normal in $J$, it is enough to show that $x^g$ belongs to $X$ for every $x\in H_1^\sigma\cup\ldots\cup H_n^\sigma$ and $g\in H_1\cup\ldots\cup H_n$. Suppose $x\in H_\ell^\sigma$ and $g\in H_m$, where $1\leq \ell,m\leq n$. Then Theorem \ref{theopermuta} yields that $x^g\in \big(H_\ell H_m\big)^\sigma=H_\ell^\sigma H_m^\sigma\leq X$. This shows that $X$ is normal in $J$.

Since $J^\sigma\geq H_j^\sigma$ for every $j$, it follows that~\hbox{$J^\sigma\geq X$.} It only remains to prove that~\hbox{$X\leq J^\sigma$.} To this aim, we note that the factor group $J/X$ is generated by finitely many~\hbox{$\sigma$-sub}\-normal $\sigma$-nilpotent subgroups, so Lemma 2.6 of \cite{MFsigma} gives that $J/X$ is~\hbox{$\sigma$-nil}\-po\-tent. Thus~\hbox{$J^\sigma\leq X$} and the proof is complete.
\end{proof}

\begin{cor}\label{corollaryinf}
Let $\sigma=\{\sigma_i\,:\, i\in I\}$ be a partition of $\mathbb{P}$. Let $\mathcal{H}$ be a family of $\sigma$-subnormal subgroups of a locally finite group $G$. If $J=\langle H\,:\, H\in\mathcal{H}\rangle$, then $J^\sigma=\langle H^\sigma\,:\, H\in\mathcal{H}\rangle$.
\end{cor}
\begin{proof}
Let $X=\langle H^\sigma\,:\, H\in\mathcal{H}\rangle$. Since $J^\sigma\geq H^\sigma$ for all $H\in\mathcal{H}$, we have $J^\sigma\geq X$. If~\hbox{$g\in J$} and $x\in X$, then Corollary \ref{corfinitelymanysub} yields that there are subgroups $H_1,\ldots, H_n\in \mathcal{H}$ such that $g\in K$ and $x\in K^\sigma$, where $K=\langle H_1,\ldots,H_n\rangle$; thus $x^g\in K^\sigma=H_1^\sigma\ldots H_n^\sigma\leq X$. The arbitrariness of $x$ and $g$ shows that $X$ is normal in $J$. Finally, $J/X$ is locally~\hbox{$\sigma$-hyper}\-central (again by Corollary \ref{corfinitelymanysub}) and so even $\sigma$-hypercentral. Then $J^\sigma\leq X$ and consequently $X=J^\sigma$. The statement is proved.
\end{proof}

\begin{cor}\label{teoremaadolfoqualcosa}
Let $\sigma=\{\sigma_i\,:\, i\in I\}$ be a partition of $\mathbb{P}$, and let $\tau$ be any subset of $\sigma$. If $H$ and $K$ are~\hbox{$\sigma$-sub}\-normal subgroups of a locally finite group $G=\langle H,K\rangle$, then $G^\tau=\langle H^\tau,\,K^\tau\rangle$.
\end{cor}
\begin{proof}
It follows from Theorem \ref{theopermuta} that $$G^\sigma=\langle H^\sigma,K^\sigma\rangle\leq\langle H^\tau,\,K^\tau\rangle\cap G^\tau.$$ Moreover, $(HG^\sigma/G^\sigma)^\tau=H^\tau G^\sigma/G^\sigma$ and $(KG^\sigma/G^\sigma)^\tau=K^\tau G^\sigma/G^\sigma$, so we may assume~$G$ is $\sigma$-hypercentral. Write $G=X\times Y$, where $Y$ is the largest $\tau$-hypercentral subgroup of $G$; in particular, $X=G^\tau$. Since $X=\langle X\cap H,\, X\cap K\rangle$, we have $$\langle (X\cap H)^\tau,\, (X\cap K)^\tau\rangle=X.$$ Moreover $(Y\cap H)^\tau=\{1\}=(Y\cap K)^\tau$ and so $\langle H^\tau,\, K^\tau\rangle=G^\tau$.
\end{proof}

\medskip

In Corollary \ref{teoremaadolfoqualcosa}, we cannot have always $G^\tau=H^\tau K^\tau$. In fact, in such circumstances, every pair of $\sigma$-subnormal subgroups $H$ and $K$ of a locally finite group $G$ would always contain a pair of $\sigma$-embedded subgroups $L\leq H$ and $M\leq K$ such that $$X=LM=ML\trianglelefteq\trianglelefteq G\quad{\rm and}\quad X\trianglelefteq J$$ (see the proof of~Co\-rol\-la\-ry~\ref{corollaryperm2}); it would follow that $X$ is $\sigma$-embedded in $J$ (by Lem\-ma~2.6 of \cite{MFsigma}) and so~$J$ is~\hbox{$\sigma$-sub}\-normal in $G$ by Theorem 3.6 of \cite{MFsigma}. However, the examples provided in \cite{MFsigma} show that this is not always the case.

\subsection*{The classes $\mathfrak{S},\mathfrak{S}^\infty,\mathfrak{S}_\sigma,\mathfrak{S}_\sigma^\infty,\mathfrak{M}$ and groups of finite rank}

Following \cite{lennox}, we denote by $\mathfrak{S}$ (resp. $\mathfrak{S}^{\infty}$) the class of all groups in which the join of any pair (resp. any family) of subnormal subgroups is subnormal. Let $\sigma=\{\sigma_i\,:\, i\in I\}$ be a partition of $\mathbb{P}$. We define $\mathfrak{S}_\sigma$ (resp. $\mathfrak{S}^\infty_\sigma$) as the class of all locally finite groups in which the join of any pair (resp. any family) of $\sigma$-subnormal subgroups is~\hbox{$\sigma$-sub}\-normal. Furthermore, following \cite{lennox}, we denote by $\mathfrak{M}$ the class of all groups having a finite series whose factors are in $\operatorname{Min}\textnormal{-}\operatorname{sn}$ or in $\operatorname{Max}\textnormal{-}\operatorname{sn}$, \hbox{i.e.} satisfy the minimal or the maximal condition on subnormal subgroups. Recall also that the {\it subnormal socle} of a group $G$ is the subgroup $S$ generated by all the minimal subnormal subgroups of $G$. It is well known that $S=S_0\times S_1$, where $S_0$ is the direct product of all non-abelian minimal subnormal subgroups, and $S_1$ is a periodic Baer group (see \cite{Rob72}, p.181); moreover, if $G$ satisfies $\operatorname{Min}\textnormal{-}\operatorname{sn}$, then $S$ is the direct product of a finite nilpotent group and a finite number of non-abelian simple groups (see \cite{Rob72}, Lemma 5.47). 

Note that a locally finite, soluble group in $\mathfrak{M}$ is \v Cernikov, so has finite rank (a group is said to have {\it finite rank} $r$ if any of its finitely generated subgroups can be generated by $r$ elements). With respect to the class $\mathfrak{M}$, we note the following result.

\begin{lem}\label{simple}
Let $\sigma=\{\sigma_i\,:\, i\in I\}$ be a partition of $\mathbb{P}$, and let $G$ be a hyper-$\sigma$ group in $\mathfrak{M}$. Then $G$ is $\sigma$-soluble.
\end{lem}
\begin{proof}
In order to prove the result, it is enough to consider separately the case in which $G$ satisfies the maximal condition on subnormal subgroups and that in which~$G$ satisfies the minimal condition on subnormal subgroups. Moreover, if $G$ satisfies the maximal condition on subnormal subgroups, then the result is pretty clear, so we assume $G\in\operatorname{Min}\textnormal{-}\operatorname{sn}$. Let $C$ be the intersection of the centralizers in $G$ of the finite normal sections of $G$; in particular, $|G:C|$ is finite, and we may replace $G$ by $C$, thus assuming that every finite normal section of $G$ is central.

Now, let $\{S_\alpha\}$ be the normal series defined as follows: $S_0=\{1\}$; for each ordinal~$\alpha$, $S_{\alpha+1}/S_\alpha$ is the subgroup generated  by all the non-abelian minimal subnormal subgroups of $G/S_\alpha$ (it is actually the direct product of these subgroups); for every limit ordinal $\lambda$, $S_\lambda$ is the union of the $S_\alpha$ with $\alpha<\lambda$. Let $S=\langle S_\alpha\,:\, \alpha<\lambda\rangle$. Let $i\in I$ be such that $2\in\sigma_i$. Since every non-abelian locally finite simple group contains a $2$-element, it follows that $S$ is contained in $O_{\sigma_i}(G)$ (this actually also depends on the fact that $G$ is a hyper-$\sigma$ group). Thus, assuming $O_{\sigma_i}(G)=\{1\}$, we also have $S=\{1\}$.

Let $\{T_\beta\}$ be the upper subnormal socles series of $G$. If every factor of this series is abelian, then it is also finite and so central in $G$. In this case, $G$ is hypercentral, so even~\v Cernikov (by the minimal condition on subnormal subgroups) and the statement is proved. Otherwise there is some factor that is not finite; let $\gamma$ be the smallest ordinal such that $T_{\gamma+1}/T_\gamma$ is not finite. Now, $T_\gamma$ is hypercentral and so \v Cernikov (by the minimal condition on subnormal subgroups). But the finite residual of a~\v Cernikov group is covered by finite characteristic subgroups, so $T_\gamma\leq\zeta_{n}(G)$ for some positive integer~$n$. Write $T_{\gamma+1}/T_\gamma$ as the direct product of a finite subgroup $L/T_\gamma$ and finitely many non-abelian simple subgroups. Since $L\leq\zeta_{n+1}(G)$ and $\pi\big(T_{\gamma+1}/L\big)\subseteq\sigma_i$, it follows from a classical theorem of Baer that $\gamma_{n+2}\big(T_{\gamma+1}\big)\leq O_{\sigma_i}(G)=\{1\}$ (see~\cite{Rob72}, Corollary~2 to Theorem 4.21), so $T_{\gamma+1}$ is nilpotent and we get a contradiction.
\end{proof}

\medskip

It should be remarked that a similar result does not hold if we replace the class $\mathfrak{M}$ by the class of groups of finite rank (see for instance the example just before The\-o\-rem~3.10 of \cite{MFsigma}).

\subsection*{Permutability}

In what follows, if $H$ and $K$ are subgroups of a group, we denote by $P_H(K)$ the {\it permutizer} of $K$ in $H$, \hbox{i.e.} the largest subgroup of $H$ permuting with $K$. This subgroup is very relevant in the study of the join of subnormal subgroups and it turns for instance that if $H$ and $K$ are subnormal, then also $P_H(K)$ is subnormal (see The\-o\-rem~1.6.9 of~\cite{lennox}). Here, we need a corresponding result dealing with $\sigma$-subnormal subgroups.

\begin{theo}\label{permutizer}
Let $\sigma=\{\sigma_i\,:\, i\in I\}$ be a partition of $\mathbb{P}$. Let $H$ and $K$ be $\sigma$-subnormal subgroups of a locally finite group $G$. Then $P_H(K)$ is $\sigma$-subnormal in $G$.
\end{theo}
\begin{proof}
Let $P=P_H(K)$. It follows from Theorem \ref{theopermuta} that $H^\sigma\leq P$. Let $$H=U_0\trianglelefteq_\sigma U_1\trianglelefteq_\sigma\ldots\trianglelefteq_\sigma U_m=G$$ be a $\sigma$-series connecting $H$ to $G$ and let $\sigma_1,\ldots,\sigma_\ell$ be the elements of $\sigma$ involved in this $\sigma$-series. Write $$H/H^\sigma=H_1/H^\sigma\times H_2/H^\sigma,$$ where $H_1/H^\sigma$ is the largest $\{\sigma_1,\ldots,\sigma_\ell\}$-nilpotent subgroup of $H/H^\sigma$. It follows from Lem\-ma~\ref{alternative} (and the remark before the statement) that $H$ contains a subnormal subgroup $S$ of $G$ such that $S\trianglelefteq H$ and $H/S$ is $\{\sigma_1,\ldots,\sigma_\ell\}$-nilpotent, so in particular $H_2\leq S$. It follows that $H_2$ is subnormal in $G$. 

Now, $$P/H^\sigma=P_1/H^\sigma\times P_2/H^\sigma,$$ where $P_1\leq H_1$ and $P_2\leq H_2$. It is clear that $P_1\trianglelefteq\trianglelefteq_\sigma H_1\trianglelefteq H\trianglelefteq\trianglelefteq_\sigma G$, so if we can prove that also $P_2$ is $\sigma$-subnormal in $G$, then Theorem 3.15 of \cite{MFsigma} yields that $P=P_1P_2$ is $\sigma$-subnormal in $G$. 

It follows from Theorem \ref{theopermuta} that $J^\sigma=H^\sigma K^\sigma$, where $J=\langle H,K\rangle$. Since $PKJ^\sigma=KPJ^\sigma$ and $J/J^\sigma$ is $\sigma$-hypercentral, we have that $$P_2K=P_2H^\sigma K^\sigma K=P_2KJ^\sigma=KP_2J^\sigma=KK^\sigma H^\sigma P_2=KP_2.$$ Therefore $P_2=P_{H_2}(K)$ and it turns out that we only need to prove the statement whenever $H$ is subnormal.

We proceed by induction on the defect $m$ of $H$. The statement is clear if~\hbox{$m\leq 1$,} so we may assume $m>1$. Let $L=H^G$ and $M=P^K\cap K$. Then $M$ is~\hbox{$\sigma$-sub}\-normal in~$G$ and is contained in $L$, so we may apply induction on $m$ to get that $Q=P_H(M)$ is~\hbox{$\sigma$-sub}\-normal in $L$ and so in $G$; since $P^K=PM$, we have that $P\leq Q$. Let $R=QM$ and notice that $P^K=PM\leq QM=R$. Since $R$ is $\sigma$-subnormal in~$G$ by Theorem 3.15 of \cite{MFsigma}, there is a finite $\sigma$-series $$R=R_0\trianglelefteq_\sigma R_1\trianglelefteq_\sigma\ldots\trianglelefteq_\sigma R_n=G$$ connecting $R$ to~$G$. It follows from Lemma 2.5 of \cite{MFsigma} that there is a $\sigma$-subnormal subgroup $X$ of~$G$ such that $X$ is $K$-invariant and $P^K\leq X\leq R$. Now, $X=(Q\cap X)M$, so~$XK=(Q\cap X)K$ and hence $Q\cap X\leq P$. But $P\leq Q\cap X$, and so $P=Q\cap X$ is~\hbox{$\sigma$-sub}\-normal in $G$.
\end{proof}

\medskip

We end this section with the following example. Let $p$ and $q$ be distinct prime numbers, and let $\sigma=\big\{\{p,q\},\mathbb{P}\setminus\{p,q\}\big\}$. Let $G=H\wr K$, where $H$ and $K$ are cyclic groups of orders $p$ and $q$, respectively. Then $HK\neq KH$ but $H^\sigma K^\sigma=K^\sigma H^\sigma$, so we do not have an analogue of Theorem 3.6 for the permutability of two $\sigma$-subnormal subgroups.

\section{The classes $\mathfrak{S}_\sigma^\infty$ and $\mathfrak{S}_{\operatorname{sn}\rightarrow\sigma}^\infty$}\label{sec3}

In this section we deal with the class $\mathfrak{S}_\sigma^\infty$ (resp. $\mathfrak{S}_{\operatorname{sn}\rightarrow\sigma}^\infty$) of all locally finite groups in which the join of arbitrarily many $\sigma$-subnormal subgroups (resp. subnormal subgroups) is $\sigma$-subnormal. A couple of preliminary remarks are in order. The example just before Theorem 3.10 of \cite{MFsigma} shows that $\mathfrak{S}_\sigma^\infty$ is strictly contained in $\mathfrak{S}_{\operatorname{sn}\rightarrow\sigma}^\infty$. Moreover, $\mathfrak{S}_\sigma^\infty$ is incomparable with the class $\mathfrak{S}^\infty$: if $G$ is a locally finite, metabelian primary group in which the join of arbitrarily many subnormal subgroups is not always subnormal, then $G\in\mathfrak{S}_\sigma^\infty\setminus\mathfrak{S}^\infty$; if $G$ is the example just before Theorem 3.10 of \cite{MFsigma}, then $G$ belongs to $\mathfrak{S}^\infty\setminus\mathfrak{S}_\sigma^\infty$. A diagram showing the relationship among the previous classes of groups can be found at the beginning of Section \ref{sec4}.

\medskip

In order to prove that $\mathfrak{M}\mathfrak{S}_\sigma^\infty\mathfrak{M}=\mathfrak{S}_\sigma^\infty$, we need to characterize the class $\mathfrak{S}^\infty_\sigma$ in terms of unions of ascending chains of $\sigma$-subnormal subgroups (see Theorem \ref{theomaximal}). This is accomplished through the following lemma: it is an analog of \cite{lennox}, Lemma 1.3.5, but, as you can see, although the proof in the subnormal case is almost trivial, in case of~\hbox{$\sigma$-sub}\-normal subgroups is not easy.

\begin{lem}\label{maximal}
Let $\sigma=\{\sigma_i\,:\, i\in I\}$ be a partition of $\mathbb{P}$, and let $S$ be a subgroup of the locally finite group $G$. If $H$ is a maximal member of the set of $\sigma$-subnormal subgroups of $G$ lying in~$S$, then $H\trianglelefteq S$.
\end{lem}
\begin{proof}
Let $$H=H_0\trianglelefteq_\sigma H_1\trianglelefteq_\sigma H_2\trianglelefteq_\sigma\ldots\trianglelefteq_\sigma H_n=S$$ be a $\sigma$-series connecting $H$ and $S$ of smallest possible length $n$. Assume the result is false.

If $n=1$, then $H\trianglelefteq_{\sigma_i} S$, for some $i\in I$. Let $x\in S$ be such that $H_0^x\not\leq H_0$. Since $H_0$ and~$H_0^x$ are $\sigma$-subnormal subgroups of $G$, and $H_0$ is $\sigma$-normal in $J=\langle H_0,H_0^x\rangle$ (because $S\geq J$), it follows from Lemma \ref{cor3.14} that $J$ is $\sigma$-subnormal in $G$. The maximality of $H_0$ gives the contradiction $H_0=J\geq H_0^x$. Therefore $n\geq 2$.

Now, assume $H$ is not normal in $H_1$, and choose $x\in H_1$ such that $H^x\not\leq H$.  An argument similar to that of the previous paragraph yields that $\langle H,H^x\rangle$ is a $\sigma$-subnormal subgroup of $G$ which is contained in $H_1\leq S$. Again the maximality of $H$ gives a contradiction. Therefore $H\trianglelefteq H_1$.

We claim that $H_1$ is not normal in $H_2$. Suppose $H_1$ is normal in $H_2$, and notice that~$H$ cannot be normal in $H_2$ (by minimality of $n$). Choose $x\in H_2$ such that $H^x\not\leq H$. Then $J=\langle H,H^x\rangle$ is contained in $H_1$, and hence $J=H^xH$. It follows from Theorem 3.15 of~\cite{MFsigma} that $J$ is $\sigma$-subnormal in $G$, and again the contradiction $H=J\geq H^x$. Therefore $H_1\trianglelefteq_{\sigma_i} H_2$ for some $i\in I$.

Let $L=(H_1)_{H_2}$, so $H_2/L$ is a $\sigma_i$-group. Notice now that $H\cap L$ is a maximal~\hbox{$\sigma$-sub}\-normal subgroup of $G$ contained in $L$. In fact, $H\cap L$ is normal in $H$, and so is~\hbox{$\sigma$-sub}\-normal in $G$. If $M$ is a $\sigma$-subnormal subgroup of $G$ such that $$H\cap L\leq M\leq L,$$ then $H$ is normal in $\langle H, M\rangle\leq\langle H,L\rangle\leq H_1$, and hence $HM$ is $\sigma$-subnormal in $G$ by~Lem\-ma~\ref{cor3.14}. Thus, $M\leq H$ and so $H\cap L=M$. Being $$H\cap L\trianglelefteq L\trianglelefteq H_2,$$ an argument similar to the one employed in the fourth paragraph of the proof yields that $H\cap L$ is normal in $H_2$ ($H\cap L$ must contains all its conjugates in $H_2$).

Finally, $\pi(H/H\cap L)\subseteq\sigma_i$, so Lemma 2.6 of \cite{MFsigma} yields that $H/H\cap L$ is contained in a~\hbox{$\sigma_i$-sub}\-group $Q/H\cap L$ that is a normal subgroup of $H_2/H\cap L$. Therefore $$H\trianglelefteq_{\sigma_i} Q\trianglelefteq H_2,$$ contradicting what we showed above (here, $Q$ takes the place of $H_1$). The statement is proved.
\end{proof}

\begin{cor}\label{cormaximal}
Let $\sigma=\{\sigma_i\,:\, i\in I\}$ be a partition of $\mathbb{P}$. Let $\{H_\lambda\,:\, \lambda\in\Lambda\}$ be a set of~\hbox{$\sigma$-sub}\-normal subgroups of a locally finite group $G$ and let $J$ be their join. Then $J$ is $\sigma$-subnormal in $G$ if and only if the set of all $\sigma$-subnormal subgroups of $G$ lying in $J$ contains a maximal member.
\end{cor}
\begin{proof}
We only need to prove that the condition is sufficient. Let $M$ be a~\hbox{$\sigma$-sub}\-normal subgroup of $G$ maximal with respect to lying in $J$. By Lemma \ref{maximal}, $M$ is normal in~$J$, so, for every $\lambda\in\Lambda$, $H_\lambda M$ is $\sigma$-subnormal in $G$ by Lemma \ref{cor3.14}. It follows that $H_\lambda\leq M$ for all $\lambda\in\Lambda$ and hence $J=M$ is $\sigma$-subnormal in $G$.
\end{proof}

\medskip

We are ready to prove the promised characterization of the class $\mathfrak{S}^\infty_\sigma$ (see also~The\-o\-rem~2.4.4 of \cite{lennox}).

\begin{theo}\label{theomaximal}
Let $\sigma=\{\sigma_i\,:\, i\in I\}$ be a partition of $\mathbb{P}$, and let $G$ be a locally finite group. Then $G\in\mathfrak{S}^\infty_\sigma$ if and only if the union of each ascending chain of $\sigma$-subnormal subgroups of~$G$ is $\sigma$-subnormal in $G$.
\end{theo}
\begin{proof}
Of course, we only need to prove the sufficiency of the condition. Let $J$ be a join of $\sigma$-subnormal subgroups of $G$. By Zorn's lemma, the set of $\sigma$-subnormal subgroups of $G$ lying in $J$ contains a maximal member. Now, Corollary \ref{cormaximal} yields that~$J$ is $\sigma$-subnormal in $G$. Therefore $G$ belongs to $\mathfrak{S}^\infty_\sigma$.
\end{proof}

\medskip

We now separately deal with the class $\mathfrak{S}_\sigma^\infty\mathfrak{M}$ and the class $\mathfrak{M}\mathfrak{S}_\sigma^\infty$. We start considering proving that $\mathfrak{S}_\sigma^\infty\mathfrak{M}=\mathfrak{S}_\sigma^\infty$. Here, we preliminary need to prove a corresponding result for the class $\mathfrak{S}^\infty_{\operatorname{sn}\rightarrow\sigma}$.

\begin{theo}\label{theomaximalsnsigma}
Let $\sigma=\{\sigma_i\,:\, i\in I\}$ be a partition of $\mathbb{P}$, and let $G$ be a locally finite group. Then $G\in\mathfrak{S}^\infty_{\operatorname{sn}\rightarrow\sigma}$ if and only if the union of each ascending chain of subnormal subgroups of~$G$ is $\sigma$-subnormal in $G$.
\end{theo}
\begin{proof}
Only the sufficiency is in question. Let $\mathcal{H}$ be a family of subnormal subgroups of $G$ and let $J$ be their join. The Axiom of Choice implies that $\mathcal{H}$ can be well-ordered; so there is an ordinal number $\lambda$ such that $\mathcal{H}=\{H_\alpha\,:\, \alpha<\lambda\}$. Now, for each ordinal number $\beta\leq\lambda$, we put $X_\beta=\langle H_\alpha\,:\,\alpha<\beta\rangle$. Of course, $X_0=\{1\}$, $X_1=H_0$ is subnormal in $G$, and $X_\lambda=J$. Assume by contradiction that $J$ is not $\sigma$-subnormal in~$G$ and let $\mu$ be the smallest ordinal number such that $X_\mu$ is not~\hbox{$\sigma$-sub}\-normal in $G$. It is clearly possible to assume $\lambda=\mu>1$.

For every $\alpha<\lambda$, the $\sigma$-residual $X_\alpha^\sigma$ is a subnormal subgroup of $G$ by Lemma \ref{alternative} since~$X_\alpha$ is $\sigma$-subnormal in $G$. If $\lambda$ is successor, then $J$ is the join of the $\sigma$-subnormal subgroups $X_{\lambda-1}$ and $H_{\lambda-1}$, so $J^\sigma$ is subnormal in $G$ by~Co\-rol\-la\-ry~\ref{corollaryperm2}. If $\lambda$ is limit, the $\sigma$-residual $J^\sigma$ of $J$ is the union of the chain of subnormal sub\-groups $\{X_\alpha^\sigma\}_{\alpha<\lambda}$ (you can see this using the argument in the proof of~The\-o\-rem~\ref{theopermuta}), which means that~$J^\sigma$ is~\hbox{$\sigma$-sub}\-normal in $G$. In any case $J^\sigma$ is~\hbox{$\sigma$-sub}\-normal in $G$. By~Lem\-ma~\ref{alternative} (and the remark before its statement) we can find a characteristic sub\-group~$W$ of $J^\sigma$ that is~\hbox{$\sigma$-embedded} in $J^\sigma$ and subnormal in $G$; let $\ell$ be the subnormal defect of~$W$ in $G$. We work by induction on $\ell$ (\hbox{i.e.} on the defect of a subnormal subgroup that is~\hbox{$\sigma$-embedded} in a normal subgroup containing the $\sigma$-residual). Being the statement obvious when $\ell=0$, assume~\hbox{$\ell>0$.} Let $$W=W_0\trianglelefteq W_1\trianglelefteq\ldots\trianglelefteq W_\ell=G$$ be the normal closure series of $W$ in $G$. Since $J$ is normalizes $W$, then it also normalizes~$W_1$. Moreover, for each $\alpha<\lambda$, $H_\alpha W_1$ is subnormal in $G$; in particular,~$W_1J$ is a join of subnormal subgroups of $G$. But~$W_1$ is normal in $W_1J$ and $W_1J/W_1\simeq J/W_1\cap J$ is a quotient of $J/W$, so the induction hypothesis yields that $W_1J$ is $\sigma$-subnormal in $G$.  

Theorem 3.10 of \cite{MFsigma} shows that $J$ is $\sigma$-serial in $G$, so Lemma 2.6 of \cite{MFsigma} yields that~\hbox{$T/W_0=J^{W_1}/W_0$} contains a normal $\sigma$-soluble subgroup $M/W_0$ such that~\hbox{$T/M$} is $\sigma$-hypercentral. Let $R$ be a subnormal subgroup of $G$ that is $\sigma$-embedded in $T$ and normal in $T$. Replacing $R$ by $W_0R$, we may assume~$R\geq W_0$. We claim that $(M\cap R)/W_0$ is a normal $\sigma$-soluble subgroup of $T/W_0$  such that $R/(M\cap R)\in\mathfrak{S}^\infty$. Since $R/(M\cap R)$ is $\sigma$-hypercentral, we can write $$R/(M\cap R)=R_1/(M\cap R)\times R_2/(M\cap R)\times\ldots \times R_m/(M\cap R)\times\ldots$$ as the direct product of its $\sigma$-components. Assume there is a strictly increasing series of natural numbers $j_1<j_2<\ldots$ such that, for each $n$, there is a subnormal subgroup $S_{j_n}/(M\cap R)$ of  $R_{j_n}/(M\cap R)$ of defect at least $n$. In this case, $\langle S_{j_1},S_{j_2},\ldots\rangle$ is union of subnormal subgroups of $G$ whose join is not $\sigma$-subnormal, a contradiction. It is therefore possible to assume that there is an integer $m$ such that each subnormal subgroup of  $R/(M\cap R)$ has defect at most $m$. Now, the union of any chain of subnormal subgroups of $R/(M\cap R)$ is easily seen to be subnormal, so~\hbox{$R/(M\cap R)$} is in $\mathfrak{S}^\infty$ by~The\-o\-rem~2.4.4 of \cite{lennox}. The claim is proved.

Now, it follows from Corollary 3.5 of \cite{MFsigma} that $J$ is $\sigma$-subnormal in $JM$. The previous claim yields that we can write $T/M=A/M\times B/M$, where \hbox{$\pi(A/M)\cap\pi(B/M)=\emptyset$}, $A/M$ is $\sigma$-nilpotent and $B/M\in\mathfrak{S}^\infty$. Thus, $A\cap JM$ is $\sigma$-subnormal in $T/M$ by~The\-o\-rem~3.6 of \cite{MFsigma} and, being a join of subnormal subgroups, $B\cap JM$ is subnormal in $B$ and so in $T$. Finally, $JM=(A\cap JM)(B\cap JM)$ is~\hbox{$\sigma$-sub}\-normal in $T$ and so in $G$ by~The\-o\-rem~3.15 of \cite{MFsigma}. This contradiction completes the proof of the result.
\end{proof}

\begin{theo}\label{sinftymaxsnpre}
Let $\sigma=\{\sigma_i\,:\, i\in I\}$ be a partition of $\mathbb{P}$. Let $G$ be a locally finite group having a nor\-mal~$\mathfrak{S}^\infty_{\operatorname{sn}\rightarrow\sigma}$-subgroup $N$ such that $G/N\in\mathfrak{M}$. Then $G\in\mathfrak{S}^\infty_{\operatorname{sn}\rightarrow\sigma}$.
\end{theo}
\begin{proof}
By Theorem \ref{theomaximalsnsigma}, we only need to prove that the union of an ascending chain of subnormal subgroups is $\sigma$-subnormal. Let $$V_1\leq V_2\leq\ldots V_\alpha\leq V_{\alpha+1}\leq\ldots\quad\quad (\alpha<\mu)$$ be  an ascending chain of subnormal subgroups and let $V$ be their union. Since $VN/N$ is subnormal in $G/N$ (see Theorem 2.6.1 of \cite{lennox}), and so $VN$ is subnormal in $G$, we may assume $G=VN$. Moreover, $V\cap N$ is the union of an ascending chain of subnormal subgroups of $N$, so the hypotheses yield that $V\cap N$ is $\sigma$-subnormal (in $N$ and so) in~$G$. Let $$V\cap N=U_0\trianglelefteq_\sigma U_1\trianglelefteq_\sigma\ldots\trianglelefteq_\sigma U_e=N$$ be a $V$-invariant $\sigma$-series connecting $V\cap N$ to $N$ (whose existence is provided by~Lem\-ma~2.5 of \cite{MFsigma}). 

If $U_0$ is normal in $U_1$, then $U_0$ is normal in $VU_1$, and 
$VU_1/U_0=V/U_0\ltimes U_1/U_0$. Let $L=C_V(U_1/U_0)$; in particular, $L\geq U_0$, $L\trianglelefteq VU_1$ and $VU_1/L= V/L\ltimes U_1L/L$. For each $\alpha<\mu$,~$V_\alpha L$ is a subnormal subgroup of $VU_1$, so $V_\alpha L/L$ stabilizes a finite chain of~$U_1\simeq U_1L/L$ and is therefore nilpotent by a well known theorem of Philip~Hall. Then~$V/L$ is locally nilpotent. But $V/V\cap N$ is in $\mathfrak{M}$, so $V/L$ is the direct product of finitely many primary components, and hence (since $V/L$ is the union of subnormal subgroups of~$VU_1/L$) the normal closure~$W/L$ of $V/L$ in $VU_1/L$ is still the direct product of finitely many primary components. This means that $W/L$ is $\sigma$-nilpotent and so that $V$ is $\sigma$-subnormal in $VU_1$.

If $U_0$ is $\sigma_i$-normal in $U_1$ for some $i\in I$, then $U_0$ contains a $VU_1$-invariant subgroup~$Q$ such that $U_1/Q$ is a $\sigma_i$-group. On the other hand, $V$ is $\sigma$-serial in $VU_1$ (see~The\-o\-rem~3.1 of \cite{MFsigma}), so Corollary 3.5 of \cite{MFsigma} shows that $V/Q\trianglelefteq_{\sigma_i} VU_1/Q$. Thus, also in this case, $V$ is $\sigma$-subnormal in $VU_1$.

Repeating the above argument, we see that $VU_i$ is $\sigma$-subnormal in $VU_{i+1}$ for each $0\leq i<e$, so $$V=VU_0\trianglelefteq\trianglelefteq_{\sigma} VU_1\trianglelefteq\trianglelefteq_{\sigma}\ldots\trianglelefteq\trianglelefteq_{\sigma} VU_e=VN=G$$ and the statement is proved.
\end{proof}

\begin{theo}\label{sinftymaxsn}
Let $\sigma=\{\sigma_i\,:\, i\in I\}$ be a partition of $\mathbb{P}$. Let $G$ be a locally finite group having a nor\-mal~$\mathfrak{S}^\infty_\sigma$-subgroup $N$ such that $G/N\in\mathfrak{M}$. Then $G\in\mathfrak{S}^\infty_\sigma$.
\end{theo}
\begin{proof}
By Theorem \ref{theomaximal}, we only need to show that the union of an ascending chain of $\sigma$-subnormal subgroups of $G$ is $\sigma$-subnormal. Let $\lambda$ be an ordinal number, and let $$H_1\leq H_2\leq\ldots H_\alpha\leq H_{\alpha+1}\leq\ldots\quad\quad (\alpha<\lambda)$$ be an ascending chain of $\sigma$-subnormal subgroups; let $H$ be the union of the subgroups in this chain.

Suppose first that the cofinality of $\lambda$ is strictly larger than $\omega$. Then there is a subsequence $\big\{H_{\alpha(\beta)}\big\}_{\beta<\lambda}$ of $\{H_\alpha\}_{\alpha<\lambda}$ whose union is still $H$ and for which there exist $i_0,\ldots,i_{n-1}\in I\cup\{I\}$ such that, for any $\beta<\lambda$, there is a $\sigma$-series $$H_{\alpha(\beta)}=H_{0,\beta}\trianglelefteq_{\sigma_{i_0}} H_{1,\beta}\trianglelefteq_{\sigma_{i_1}}\ldots\trianglelefteq_{\sigma_{i_{n-1}}} H_{n,\beta}=G$$ connecting $H_{\alpha(\beta)}$ to $G$ (in the above series, $\trianglelefteq_{\sigma_I}$ means $\trianglelefteq$). Now, let $i_{j_1},\ldots,i_{j_m}$ be the subsequence obtained from $i_0,\ldots,i_{n-1}$ by removing all terms that are equal to $I$.  For each $\beta<\lambda$, let $R_\beta$ be the $\big\{\sigma_{i_{j_1}},\ldots,\sigma_{i_{j_m}}\big\}$-residual of $H_{\alpha(\beta)}$. It follows from the proof of Lemma \ref{alternative} (see the remark before the statement) that $R_\beta$ is subnormal in~$G$. Moreover, if $\beta\leq\gamma<\lambda$, then $R_\beta\leq R_\gamma$. Let $$R=\bigcup_{\beta<\lambda}R_\beta,$$ so $R$ is a subgroup of $H$ which is $\sigma$-subnormal in $G$ (by Theorem \ref{sinftymaxsnpre}).

Assume now the cofinality of $\lambda$ is $\omega$. Since we are interested in the $\sigma$-subnormality of $H$, we may as well assume $\lambda=\omega$ in this case. For every integer $k$, let $S_k$ be a subnormal subgroup of $G$ that is $\sigma$-normally embedded in~$H_k$ (see Lemma 2.9 of~\cite{MFsigma}); moreover, let $$T_k=\langle S_h\,:\, h\leq k\rangle=S_1S_2\ldots S_k\quad\textnormal{ and }\quad R=\bigcup_{h<\omega} T_h=\langle S_h\,:\, h<\omega\rangle,$$ so $T_k$ is a subnormal subgroup of $G$, and $R$ is a $\sigma$-subnormal subgroup of $G$ (again by~The\-o\-rem \ref{sinftymaxsnpre}).

In both cases we find a $\sigma$-subnormal subgroup $R$ of $G$ that is contained in $H$ and containing,  for each $\alpha<\lambda$, a $\sigma$-normally embedded subgroup $S_\alpha$ of $H_\alpha$. Let $T$ be a subnormal subgroup of $G$ that is $\sigma$-normally embedded in $R$. We claim that $TN/N$ is~\hbox{$\sigma$-embed}\-ded in $HN/N$. For the time being, we assume $N=\{1\}$; so $G$ is in $\mathfrak{M}$, and we need to prove that $T$ is~\hbox{$\sigma$-embed}\-ded in $H$.

 Let $$T=T_0\trianglelefteq T_1\trianglelefteq\ldots\trianglelefteq T_l=G$$ be the normal closure series of $T$ in $G$. By Lemma 2.6 of \cite{MFsigma}, $(T_1\cap H)^{T_1}/T$ is hyper-$\sigma$, so it is even $\sigma$-soluble by Lemma \ref{simple}, because $(T_1\cap H)^{T_1}/T$ is in $\mathfrak{M}$. Thus, $T=T_0$ is $\sigma$-embedded in $T_1\cap H$. Similarly, we find that for each $0\leq u<l$, $T_u$ is $\sigma$-embedded in $(T_{u+1}\cap H)T_u$. This means that, for each $0\leq u<l$, $T_u\cap H$ is $\sigma$-embedded in $T_{u+1}\cap H$, and hence that $T$ is $\sigma$-embedded in $H$. The claim is proved.
 
\smallskip
 
Finally, we come back to the general case. Now, $TN/N$ is $\sigma$-embedded in $HN/N$, and this means that $T(H\cap N)/(H\cap N)$ is $\sigma$-embedded in $H/H\cap N$. Moreover, $T$ is subnormal in $G$ and $H\cap N$ is $\sigma$-subnormal in $G$, so Theorem 3.15 of \cite{MFsigma} gives that~\hbox{$T(H\cap N)$} is $\sigma$-subnormal in $G$, and hence~\hbox{$T(H\cap N)$} contains a $\sigma$-embedded subgroup~$Q$ that is subnormal in $G$ (see~Lem\-ma~\ref{alternative}). In particular,~$Q$ is $\sigma$-embedded in~$H$ and is subnormal in $G$. But~The\-o\-rem~3.10 of~\cite{MFsigma} yields that $H$ is $\sigma$-serial in $G$, so~The\-o\-rem~3.6 of~\cite{MFsigma} completes the proof.
\end{proof}

\medskip

As a consequence of the above result, we can answer Question 3.24 of \cite{MFsigma}.

\begin{cor}\label{question}
$\mathfrak{M}\leq\mathfrak{S}^\infty_\sigma$ for any partition $\sigma=\{\sigma_i\,:\, i\in I\}$ of~$\mathbb{P}$.
\end{cor}

\medskip

It is well known that the class of all locally finite groups $G$ having an abelian normal subgroup $N$ with $G/N\in\operatorname{Min}\textnormal{-}\operatorname{sn}$ is not contained in $\mathfrak{S}^\infty$ (consider for example the wreath product $C_p\wr C_{p^\infty}$, for any prime $p$). It is therefore interesting to remark that, as a consequence of Theorem \ref{sinftymaxsn}, we have that in any group in this class, the join of arbitrarily many subnormal subgroups is at least $\sigma$-subnormal for any partition $\sigma$ of~$\mathbb{P}$.

Furthermore, we observe that the proof of Theorem \ref{sinftymaxsn} makes it possible to prove that {\it if $\mathfrak{X}$ is class of groups which is closed with respect to forming subnormal subgroups, quotients and is such that every hyper-$\sigma$ $\mathfrak{X}$-group is $\sigma$-soluble (actually, we just need that every $\sigma$-hypercentral $\mathfrak{X}$-group is $\sigma$-nilpotent), then $\mathfrak{S}_{\operatorname{sn}\rightarrow\sigma}^\infty\cap\mathfrak{X}\leq\mathfrak{S}_\sigma^\infty$} (see also~The\-o\-rem~\ref{theofinaleee}).

\medskip

Our next result deals with the dual case of a group that is the extension of a group in $\mathfrak{M}$ by a group in $\mathfrak{S}_\sigma^\infty$.

\begin{theo}\label{sinftyminmax}
Let $\sigma=\{\sigma_i\,:\, i\in I\}$ be a partition of $\mathbb{P}$. Let $G$ be a locally finite group having a nor\-mal subgroup $N$ in $\mathfrak{M}$ such that $G/N\in\mathfrak{S}^\infty_\sigma$. Then $G$ is in $\mathfrak{S}^\infty_\sigma$.
\end{theo}
\begin{proof}
By Theorem \ref{theomaximal}, we only need to show that the union of an ascending chain of $\sigma$-subnormal subgroups of $G$ is $\sigma$-subnormal. Let $$H_1\leq H_2\leq\ldots H_\alpha\leq H_{\alpha+1}\leq\ldots\quad\quad (\alpha<\lambda)$$ be such an ascending chain of $\sigma$-subnormal subgroups and let $H$ be their union. Since~$H$ is generated by $\sigma$-subnormal subgroups, we have that $H$ is $\sigma$-serial in $G$ (see~The\-o\-rem~3.10 of \cite{MFsigma}).

Moreover, $$H\cap N=\bigcup_{\alpha<\lambda}(H_\alpha\cap N),$$ so it follows from Corollary \ref{question} that $H\cap N$ is $\sigma$-subnormal in $G$. Now, Lem\-ma~2.5 of~\cite{MFsigma} shows that there exists an $H$-invariant $\sigma$-series $$H\cap N=K_0\trianglelefteq_\sigma K_1\trianglelefteq_\sigma K_2\trianglelefteq_\sigma\ldots\trianglelefteq_\sigma K_n=N$$ connecting $H\cap N$ to $N$. We claim that $H$ is $\sigma$-subnormal in $HK_1$.

Suppose first there is $i\in I$ such that $H\cap N\trianglelefteq_{\sigma_i} K_1$. Let $$L=(H\cap N)_{K_1}=(H\cap N)_{HK_1},$$ so in particular $K_1/L$ is a normal $\sigma_i$-subgroup of $HK_1/L$, and hence Corollary 3.5 of~\cite{MFsigma} shows that $H=HL$ is $\sigma$-subnormal in $HK_1$, thus proving the claim in this case.

Suppose now that $H\cap N$ is normal in $K_1$, so $H\cap N$ is normal in $HK_1$. Then $$HK_1/(H\cap N)=H/(H\cap N)\ltimes K_1/(H\cap N).$$ Let $L/(H\cap N)=C_{H/(H\cap N)}\big(K_1/(H\cap N)\big)$. If $S/L$ is any subnormal subgroup of $HK_1/L$ that is contained in $H/L$, then $S/L$ centralizes a finite chain of $K_1L/L\simeq K_1$, so a well known theorem of~Philip~Hall yields that $S/L$ is nilpotent. Now, a combination of~Lem\-mas~2.6 and~2.9 of \cite{MFsigma} shows that every $H_\alpha L/L$ is a hyper-$\sigma$ group, so again Lem\-ma~2.6 of \cite{MFsigma} yields that $H/L$ is contained in a normal subgroup~$M/L$ of $HK_1/L$ that is \hbox{hyper-$\sigma$.} Of course, $M/L=H/L\ltimes (M\cap K_1)L/L$. Moreover, $$(M\cap K_1)L\trianglelefteq K_1L\trianglelefteq\trianglelefteq_\sigma G,$$ so $(M\cap K_1)L$ has a $\sigma$-embedded subgroup $T$ that is subnormal in $G$; but then $TL/L$ belongs to $\mathfrak{M}$ and so Lemma~\ref{simple} gives that $TL/L$  is $\sigma$-soluble. Therefore \hbox{$(M\cap K_1)L/L$} is $\sigma$-soluble. Finally, Co\-rol\-la\-ry~3.5 (and Lemma 2.6) of \cite{MFsigma} yields that $H$ is $\sigma$-subnormal in $M$ and so even in $HK_1$. The claim is proved.

\smallskip

An argument similar to the one employed in the previous claims makes it possible to show that $HK_i$ is $\sigma$-subnormal in $HK_{i+1}$, for each $1\leq i<n$. Now, since $$H\trianglelefteq\trianglelefteq_\sigma HK_1\trianglelefteq\trianglelefteq_\sigma HK_2\trianglelefteq\trianglelefteq_\sigma\ldots\trianglelefteq\trianglelefteq_\sigma HK_n=HN,$$ we have that $H$ is $\sigma$-subnormal in $HN$. On the other hand, $G/N\in\mathfrak{S}^\infty_\sigma$, so $HN/N$ is~\hbox{$\sigma$-sub}\-normal in $G/N$ and hence $H$ is $\sigma$-subnormal in $G$. The statement is proved.
\end{proof}

\medskip

If in the previous proof we replace $\mathfrak{S}^\infty_\sigma$ by $\mathfrak{S}^\infty_{\operatorname{sn}\rightarrow\sigma}$, and Theorem \ref{theomaximal} by Theorem \ref{theomaximalsnsigma}, we get the following result.

\begin{theo}\label{sinftyminmax22}
Let $\sigma=\{\sigma_i\,:\, i\in I\}$ be a partition of $\mathbb{P}$. Let $G$ be a locally finite group having a nor\-mal subgroup $N$ in $\mathfrak{M}$ such that $G/N\in\mathfrak{S}^\infty_{\operatorname{sn}\rightarrow\sigma}$. Then $G$ is in $\mathfrak{S}^\infty_{\operatorname{sn}\rightarrow\sigma}$.
\end{theo}

\medskip

We remark that in the proofs of Theorems \ref{sinftymaxsn}, \ref{sinftyminmax} and \ref{sinftyminmax22}, the use of~Lem\-ma~\ref{simple} is not an essential one and could been replaced by a more accurate use of~Lem\-ma~\ref{alternative}. The reason we chose to do it in this way is that $\sigma$-embedded subgroups are handier to deal with in a proof and, besides, Lemma \ref{simple} is interesting in itself.

It should also be noticed that in the statements of Theorems \ref{sinftymaxsn}, \ref{sinftyminmax} and \ref{sinftyminmax22} we cannot replace in the hypothesis~$\mathfrak{S}^\infty_\sigma$ (or $\mathfrak{S}^\infty_{\operatorname{sn}\rightarrow\sigma}$) by $\mathfrak{S}^\infty$, since a locally finite group in~$\mathfrak{S}^\infty$ is not always contained in $\mathfrak{S}^\infty_\sigma$: see the example just before Theorem 3.10 of~\cite{MFsigma}. Moreover, the same example has rank $2$, and so it shows that the class of locally finite groups of finite rank is not in general contained in $\mathfrak{S}_\sigma^\infty$ (or in $\mathfrak{S}^\infty_{\operatorname{sn}\rightarrow\sigma}$) for suitable choices of a partition $\sigma$ of $\mathbb{P}$; in particular, we cannot replace in the previous four theorems the class $\mathfrak{M}$ by the class of all groups of finite rank.

A careful look at the above proofs shows that in Theorems  \ref{sinftymaxsnpre}, \ref{sinftymaxsn}, \ref{sinftyminmax} and \ref{sinftyminmax22}) we may replace $\mathfrak{M}$ by any class of groups $\mathfrak{X}\leq\mathfrak{S}^\infty_\sigma$ which is closed with respect to forming subnormal subgroups, quotients and is such that every hyper-$\sigma$ $\mathfrak{X}$-group is $\sigma$-soluble (actually, we just need that every $\sigma$-hypercentral $\mathfrak{X}$-group is $\sigma$-nilpotent); see also the remark before Theorem \ref{sinftyminmax}. This is for instance the case of the class of all $\sigma$-groups, but given the relevance of such a class in the context of $\sigma$-subnormality we also give a direct proof of this fact.

\begin{theo}\label{analogo}
Let $\sigma=\{\sigma_i\,:\, i\in I\}$ be a partition of $\mathbb{P}$. Let $G$ be a locally finite group having normal subgroups $M\leq N$ such that $M$ and $G/N$ are $\sigma$-groups, while $N/M$ is in $\mathfrak{S}_\sigma^\infty$. Then $G\in\mathfrak{S}_\sigma^\infty$.
\end{theo}
\begin{proof}
By Theorem \ref{theomaximal}, we only need to show that the union of an ascending chain of $\sigma$-subnormal subgroups of $G$ is $\sigma$-subnormal. Let $$H_1\leq H_2\leq\ldots H_\alpha\leq H_{\alpha+1}\leq\ldots\quad\quad (\alpha<\lambda)$$ be such an ascending chain of $\sigma$-subnormal subgroups and let $H$ be their union. Theorem 3.10 of \cite{MFsigma} shows that $H$ is $\sigma$-serial in $G$, so Corollary 3.5 of \cite{MFsigma} gives that $J$ is~\hbox{$\sigma$-sub}\-normal in $JM$. It is therefore possible to assume $M=\{1\}$. In this case, $H\cap N$ is $\sigma$-subnormal in $G$, being a union of $\sigma$-subnormal subgroups of $N$. Let $$H\cap N=U_0\trianglelefteq_\sigma U_1\trianglelefteq_\sigma\ldots\trianglelefteq_\sigma U_\ell=N$$ be an $H$-invariant $\sigma$-series connecting $H\cap N$ to $N$ (the existence of this $\sigma$-series is provided by Lemma 2.5 of \cite{MFsigma}). Let $i\in I$ be such that $G/N$ is a $\sigma_i$-group.

If $U_0$ is normal in $U_1$, then $U_0$ is normal in $HU_1$. Since $H/U_0$ is a union of~\hbox{$\sigma$-sub}\-normal $\sigma_i$-subgroups of $HU_1/U_0$, it follows from Lemma 2.6 of \cite{MFsigma} that $H/U_0$ is contained in a normal $\sigma_i$-subgroup of $HU_1/H_0$, and consequently $H$ is $\sigma$-subnormal in $HU_1$. 

If $U_0$ is $\sigma_j$-normal in $U_1$, for some $j\in I$, then $U_0$ contains an $HU_1$-invariant subgroup~$Q$ such that $U_1/Q$ is a $\sigma_j$-group. It follows from Corollary 3.5 of \cite{MFsigma} that $H/Q$ is $\sigma_j$-normal in $HU_1/Q$, so again $H$ is $\sigma$-subnormal in $HU_1$. 

Repeating the above argument, we see that $HU_i\trianglelefteq\trianglelefteq_\sigma HU_{i+1}$ for each $0\leq i<\ell$, so $$H=HU_0\trianglelefteq\trianglelefteq_\sigma HU_1\trianglelefteq\trianglelefteq_\sigma\ldots\trianglelefteq\trianglelefteq_\sigma HU_\ell=HN\trianglelefteq_{\sigma_i}G$$ and the statement is proved.
\end{proof}

\medskip

Essentially the same proof of Theorem \ref{analogo} makes it possible to prove the following result.

\begin{theo}\label{analogo2}
Let $\sigma=\{\sigma_i\,:\, i\in I\}$ be a partition of $\mathbb{P}$. Let $G$ be a locally finite group having normal subgroups $M\leq N$ such that $M$ and $G/N$ are $\sigma$-groups, while $N/M$ is in~$\mathfrak{S}_{\operatorname{sn}\rightarrow\sigma}^\infty$. Then $G\in\mathfrak{S}_{\operatorname{sn}\rightarrow\sigma}^\infty$.
\end{theo}

\begin{cor}
The class $\mathfrak{S}_{\operatorname{sn}\rightarrow\sigma}^\infty$ is closed with respect to forming $\sigma$-subnormal subgroups.
\end{cor}
\begin{proof}
Let $G$ be a group in $\mathfrak{S}_{\operatorname{sn}\rightarrow\sigma}^\infty$ and let $H$ be a $\sigma$-subnormal subgroup of $G$. By~Lem\-ma~\ref{alternative}, there is a subnormal subgroup $S$ of $G$ which is $\sigma$-embedded in $H$. Then~$S$ is in $\mathfrak{S}_{\operatorname{sn}\rightarrow\sigma}^\infty$ and Theorem \ref{analogo2} yields that also $H$ is in $\mathfrak{S}_{\operatorname{sn}\rightarrow\sigma}^\infty$.
\end{proof}

\medskip

We end this section with the following open problem (see Theorem \ref{finalone}). 

\begin{quest}\label{questionrelevant}
Let $\sigma=\{\sigma_i\,:\, i\in I\}$ be a finite partition of $\mathbb{P}$. Is every locally finite group in $\mathfrak{S}_\sigma^\infty$?
\end{quest}

\medskip

One can give many partial answers to this question (for example in the soluble case using Theorem 3.10 and Corollary 3.31 of \cite{MFsigma}), but the general problem seems to be very difficult to deal with. The following result could be relevant to answer the previous question. This is actually a consequence of the remark before Theorem  \ref{sinftymaxsnpre}, but for the sake of completeness, we explicitly prove it here.

\begin{theo}\label{theofinaleee}
Let $\sigma=\{\sigma_i\,:\, i\in I\}$ be a finite partition of $\mathbb{P}$. Then $\mathfrak{S}_\sigma^\infty=\mathfrak{S}_{\operatorname{sn}\rightarrow\sigma}^\infty$.
\end{theo}
\begin{proof}
We only need to show that $\mathfrak{S}_{\operatorname{sn}\rightarrow\sigma}^\infty$ is contained in  $\mathfrak{S}_\sigma^\infty$. Let $G$ be a group in~$\mathfrak{S}_{\operatorname{sn}\rightarrow\sigma}^\infty$. By Theorem \ref{theomaximal}, we need to show that every ascending chain of~\hbox{$\sigma$-sub}\-normal subgroups of $G$ is~\hbox{$\sigma$-sub}\-normal. Let $$H_1\leq H_2\leq\ldots H_\alpha\leq H_{\alpha+1}\leq\ldots\quad\quad (\alpha<\lambda)$$ be such an ascending chain of $\sigma$-subnormal subgroups and let $H$ be their union. As we know $H^\sigma=\bigcup_{\alpha<\lambda} H_\alpha^\sigma$ is a union of subnormal subgroups of $G$ (see Lemma \ref{alternative}), which means that $H^\sigma$ is $\sigma$-subnormal in $G$. But $\sigma$ is finite, and so $H^\sigma$ is $\sigma$-embedded in $H$. By Lemma \ref{alternative} (and the remark before the statement), there is a $\sigma$-embedded subgroups of $H^\sigma$ that is subnormal in $G$. Finally, an application of \cite{MFsigma}, Theorems 3.6 and 3.10, shows that $H$ is $\sigma$-subnormal in $G$ and completes the proof.
\end{proof}

\section{$\sigma$-Subnormality criteria}\label{lastsection}

The aim of this section is to provide some useful $\sigma$-subnormality criteria for the join of two $\sigma$-subnormal subgroups of a locally finite group to be $\sigma$-subnormal. These are mostly consequences of  Theorem \ref{theopermuta} and its corollaries. The first result we prove, besides answering Question 3.17 of \cite{MFsigma}, will prove its usefulness in the next section (see for example Lemma \ref{finiteextension}).

\begin{lem}\label{question2}
Let $\sigma=\{\sigma_i\,:\, i\in I\}$ be a partition of $\mathbb{P}$. Let $H$ and $K$ be $\sigma$-subnormal subgroups of a locally finite group $G$; put $J=\langle H,K\rangle$. If there are $\sigma$-embedded subgroups~$S$ and $T$ of $H$ and $K$, respectively, whose join $U=\langle S,T\rangle$ is~\hbox{$\sigma$-sub}\-normal in $G$, then $J$ is~\hbox{$\sigma$-sub}\-normal in $G$.
\end{lem}
\begin{proof}
By Corollary \ref{corollaryperm2}, $J^\sigma$ is subnormal in $G$. By Theorem 3.15 of \cite{MFsigma},  $J^\sigma H$, $J^\sigma K$ and $J^\sigma U$ are $\sigma$-subnormal in $G$. Since $J^\sigma S$ and $J^\sigma T$ are $\sigma$-embedded in $J^\sigma H$ and~$J^\sigma K$, respectively, we may therefore replace $H$, $K$, $S$ and $T$ by $J^\sigma H$, $J^\sigma K$, $J^\sigma S$ and $J^\sigma T$, respectively; in particular, $H\cap K\geq S\cap T\geq J^\sigma$. Now, we can write $J/J^\sigma=X/J^\sigma\times Y/J^\sigma$, where $X/J^\sigma$ is $\sigma$-nilpotent, $Y/J^\sigma$ is $\sigma$-hypercentral, $\pi\big(X/J^\sigma\big)\cap\pi\big(Y/J^\sigma\big)=\emptyset$, $H\cap Y=S\cap Y$ and $K\cap Y=T\cap Y$. Moreover, $$\langle S\cap Y,\, T\cap Y\rangle\trianglelefteq\langle S,T\rangle=U\trianglelefteq\trianglelefteq_\sigma G,$$ so $Y=\langle H\cap Y,\, K\cap Y\rangle=\langle S\cap Y,\, T\cap Y\rangle$ is $\sigma$-subnormal in~$G$. But $X$ is $\sigma$-subnormal in~$G$ by Theorems 3.6 and 3.10 of \cite{MFsigma}, and consequently $J=XY$ is $\sigma$-subnormal in $G$ by~The\-o\-rem 3.15 of \cite{MFsigma}.
\end{proof}

\begin{theo}\label{corquestion2}
Let $\sigma=\{\sigma_i\,:\, i\in I\}$ be a partition of $\mathbb{P}$. Let $H$ and $K$ be $\sigma$-subnormal subgroups of a locally finite group $G$; put $J=\langle H,K\rangle$. If the join of any pair of subnormal subgroups of $G$ contained in $J$ is $\sigma$-subnormal in $J$, then $J$ is $\sigma$-subnormal in $G$.
\end{theo}
\begin{proof}
This follows from Lemmas \ref{alternative} and Lemma \ref{question2}.
\end{proof}

\begin{cor}
Every locally finite group in $\mathfrak{S}$ is also contained in $\mathfrak{S}_\sigma$ for any partition $\sigma$ of~$\mathbb{P}$.
\end{cor}

\medskip

As we remarked in the introduction, the join of two subnormal subgroups $H$ and $K$ is subnormal provided that $[H,K]$ satisfies the maximal condition. Our second main result generalizes this fact to $\sigma$-subnormal subgroups.

\begin{theo}\label{corderivato}
Let $\sigma=\{\sigma_i\,:\, i\in I\}$ be a partition of $\mathbb{P}$. Let $H$ and $K$ be $\sigma$-subnormal subgroups of a locally finite group $G$. If $[H,K]\in\mathfrak{M}$, then $J=\langle H,K\rangle$ is $\sigma$-subnormal in $G$.
\end{theo}
\begin{proof}
By Corollary \ref{corollaryperm2}, we have that $J^\sigma$ is subnormal in $G$. Let $$J^\sigma=U_0\trianglelefteq U_1\trianglelefteq\ldots\trianglelefteq U_n=G$$ be the normal closure series of $J^\sigma$ in $G$ and choose $0\leq\ell<n$. We claim that~$JU_\ell/U_\ell$ is~\hbox{$\sigma$-sub}\-normal in $JU_{\ell+1}/U_\ell$. To this aim, we may assume $U_\ell=\{1\}$, so $J$ is~\hbox{$\sigma$-hyper}\-central. Since $[H,K]$ is in $\mathfrak{M}$, it follows from Lemma \ref{simple} that $[H,K]$ is $\sigma$-nilpotent, and so also $[H,K]^{JU_{\ell+1}}=[H,K]^{U_{\ell+1}}$ is $\sigma$-nilpotent by Lemma 2.6 of \cite{MFsigma}, since $J$ (and consequently~$[H,K]$) is $\sigma$-serial in $JU_{\ell+1}$ by The\-o\-rem~3.10 of \cite{MFsigma}. Application of The\-o\-rem~3.6 of \cite{MFsigma} shows that $J\trianglelefteq\trianglelefteq_\sigma J[H,K]^{U_{\ell+1}}$. It therefore possible to assume \hbox{$[H,K]=\{1\}$,} so~\hbox{$HK=KH$} and the claim follows from Theorem 3.15 of \cite{MFsigma}.

Finally, we combine all the $\sigma$-series we obtained so far to get $$J=JU_0\trianglelefteq\trianglelefteq_\sigma JU_1\trianglelefteq\trianglelefteq_\sigma\ldots\trianglelefteq\trianglelefteq_\sigma JU_n=G.$$ The statement is proved.
\end{proof}

\medskip

It is clear that in the statement of Theorem \ref{corderivato} (essentially with the same proof) we can replace $\mathfrak{M}$ by any quotient closed class of groups $\mathfrak{X}$ such that any $\sigma$-hypercentral $\mathfrak{X}$-group is~\hbox{$\sigma$-nilpotent.}

\medskip

In order to establish the $\sigma$-subnormality of a certain join $J$ of $\sigma$-subnormal subgroups of a locally finite group $G$, it can be useful to have some information about those subgroups of $J$ that are $\sigma$-subnormal in $G$. Our next results go in this direction and generalize well known theorems of Roseblade (see \cite{lennox}, Theorems 1.6.10 and 1.6.11).

\begin{theo}\label{theocontenuto}
Let $\sigma=\{\sigma_i\,:\, i\in I\}$ be a partition of $\mathbb{P}$. Let $\mathcal{L}$ be a finite family of~\hbox{$\sigma$-sub}\-normal subgroups of a locally finite group $G$. If $F$ is a subset of $J=\langle L\,:\, L\in\mathcal{L}\rangle$ such that~$\pi(F)$ is covered by finitely many elements of $\sigma$, then there is a subnormal subgroup $X$ of~$G$ such that $F\leq X\leq J$.
\end{theo}
\begin{proof}
Write $\mathcal{L}=\{L_1,\ldots,L_n\}$. It follows from Lemma \ref{alternative} and Corollary \ref{corfinitelymanysub} that $J^\sigma=L_1^\sigma\ldots L_n^\sigma$ is subnormal in $G$. Moreover, $\langle F\rangle\, J^\sigma/J^\sigma$ is contained in a normal~\hbox{$\sigma$-nil}\-potent subgroup of $J/J^\sigma$, so $\langle F\rangle\,J^\sigma$ is~\hbox{$\sigma$-sub}\-normal in $J$ and consequently $\sigma$-serial in $G$, since $J$ is $\sigma$-serial in $G$ by~The\-o\-rem~3.10 of \cite{MFsigma}. Finally, it follows from Theorem 3.6 of~\cite{MFsigma} that $\langle F\rangle\,J^\sigma$ is subnormal in $G$. The statement is proved.
\end{proof}

\medskip

Of course, the above statement applies in particular whenever $F$ is finite, and in such circumstances we can actually remove the hypothesis that $\mathcal{L}$ is finite. However, we do not know if one can remove the finiteness assumption on $\mathcal{L}$ without strengthening the condition on $F$ (this problem seems to be connected with Question \ref{questionrelevant}). 

\begin{cor}
Let $\sigma=\{\sigma_i\,:\, i\in I\}$ be a partition of $\mathbb{P}$. Let $\mathcal{L}$ be a family of $\sigma$-subnormal subgroups of a locally finite group $G$. If $H$ is a subgroup of $G$ such that $N\in\mathfrak{S}_\sigma^\infty$ and $J=\langle L\,:\, L\in\mathcal{L}\rangle$, then $J\cap H$ is $\sigma$-subnormal in $H$.
\end{cor}
\begin{proof}
By Theorem \ref{theocontenuto}, for any $x\in J\cap H$, there is a $\sigma$-subnormal subgroup $H_x$ of $G$ such that $x\in H_x\leq J$. Then $H\cap H_x$ is $\sigma$-subnormal in $H$ and contains $x$. The arbitrariness of $x$ in $J\cap H$ and the fact that $H\in\mathfrak{S}_\sigma^\infty$ show that $J\cap H$ is $\sigma$-subnormal in~$H$.
\end{proof}

\medskip

Actually, the hypothesis on $H$ in the previous statement can be slightly weakened.

\begin{theo}\label{lemmajoinfacile}
Let $\sigma=\{\sigma_i\,:\, i\in I\}$ be a partition of $\mathbb{P}$. Let $G$ be a locally finite group,~$\mathcal{L}$ a finite family of $\sigma$-subnormal subgroups of $G$, and $H$ a subgroup of $G$ in $\mathfrak{S}_{\operatorname{sn}\rightarrow\sigma}^\infty$. If \hbox{$J=\langle L\,:\, L\in\mathcal{L}\rangle$,} then $H\cap J$ is $\sigma$-subnormal in $H$. 
\end{theo}
\begin{proof}
Let $S$ be the subgroup generated by all subnormal subgroups of $G$ contained in $J$ (in particular, $S$ is normal in $J$). If $x\in H\cap S$, then there is a subnormal subgroup $S_x$ of $G$ such that $x\in S_x\leq S$ (see \cite{lennox}, Theorem 1.6.14); in particular, $x$ is contained in the subnormal subgroup $H\cap S_x$ ($\leq H\cap S$) of $H$. Therefore $H\cap S$ is~\hbox{$\sigma$-sub}\-normal in $H$ by hypothesis. 

Now, $$(H\cap J)/(H\cap S)\simeq(H\cap J)S/S,$$ so $H\cap S$ is $\sigma$-embedded in $H\cap J$ by Lemma 2.6 of \cite{MFsigma} (recall also that $\mathcal{L}$ is finite).  Finally,~Lem\-ma~\ref{alternative} shows that~\hbox{$H\cap S$} contains a $\sigma$-embedded subgroup $T$ that is subnormal in $H$, so $T$ is $\sigma$-embedded in $H\cap J$ ; but $J$ is $\sigma$-serial in $G$ (see~The\-o\-rem~3.10 of~\cite{MFsigma}), so $H\cap J$ is $\sigma$-serial in $H$ and hence Theorem 3.6 of \cite{MFsigma} yields that $H\cap J$ is~\hbox{$\sigma$-sub}\-normal in~$H$.
\end{proof}


\begin{theo}\label{altroteorem}
Let $\sigma=\{\sigma_i\,:\, i\in I\}$ be a partition of $\mathbb{P}$. Let $H$ and $K$ be $\sigma$-subnormal subgroups of a locally finite group $G$. If $L\leq H$ and $M\leq K$ are such that $LM=ML$, then there is a $\sigma$-subnormal subgroup $X$ of $G$ such that $LM\leq X\leq HK$.
\end{theo}
\begin{proof}
Let $J=\langle H,K\rangle$. By Corollary \ref{corollaryperm2}, $J^\sigma$ is subnormal in $G$. On the other hand,~Lem\-ma~\ref{alternative} shows that $H^\sigma$ and $K^\sigma$ are subnormal in $G$, so $J^\sigma H$ and $J^\sigma K$ are $\sigma$-subnormal subgroups of $G$ (see \cite{MFsigma}, Theorem 3.15). Moreover, $$J^\sigma H J^\sigma K=HJ^\sigma K=HH^\sigma K^\sigma K=HK$$ by Theorem \ref{theopermuta}.  It is therefore safe to assume $H\cap K\geq J^\sigma$ and actually even that $L\cap M\geq J^\sigma$.

As in Theorem \ref{permutizer}, we can write $$H/H^\sigma=H_1/H^\sigma\times H_2/H^\sigma\quad\textnormal{and}\quad K/K^\sigma=K_1/K^\sigma\times K_2/K^\sigma$$ where $H_1/H^\sigma$ and $K_1/K^\sigma$ are $\sigma$-nilpotent, $H_2$ and $K_2$ are subnormal in $G$, and $$
\begin{array}{c}
\pi\big(H_1/H^\sigma\big)\cap\pi\big(H_2/H^\sigma\big)=\pi\big(K_1/K^\sigma\big)\cap\pi\big(K_2/K^\sigma\big)\\[0.2cm]
=\pi\big(H_1/H^\sigma\big)\cap\pi\big(K_2/K^\sigma\big)=\pi\big(K_1/K^\sigma\big)\cap\pi\big(H_2/H^\sigma\big)=\emptyset.
\end{array}
$$ Now, $\pi\big(H_1/H^\sigma\big)\cup\pi\big(K_1/K^\sigma\big)$ is covered by finitely many elements $\sigma_1,\ldots,\sigma_\ell$ of $\sigma$. Let $\pi_1=\sigma_1\cup\ldots\cup\sigma_\ell$. Then we can write $$J/J^\sigma=J_1/J^\sigma\times J_2/J^\sigma,$$ where $J_1/J^\sigma$ is the $\pi_1$-component of $J/J^\sigma$, so $J_1/J^\sigma$ is $\sigma$-nilpotent. It follows that $$J_1=\langle H_1,K_1,J^\sigma\rangle\quad\textnormal{ and }\quad J_2=\langle H_2,K_2,J^\sigma\rangle.$$ Write $$L/J^\sigma=L_1/J^\sigma\times L_2/J^\sigma\quad\textnormal{and}\quad M/J^\sigma=M_1/J^\sigma\times M_2/J^\sigma,$$ where $L_1\cup M_1\subseteq J_1$ and $L_2\cup M_2\subseteq J_2$.  Clearly, $M_1L_1=L_1M_1$ and $M_2L_2=L_2M_2$. 

Now, $J^\sigma$ is subnormal in $G$ and $\sigma$-embedded in $M_1L_1$, so $M_1L_1$ is $\sigma$-subnormal in~$G$ by Theorem~3.6 of \cite{MFsigma}, because $J$ (so $J_1$ and hence $M_1L_1$) is $\sigma$-serial in $G$ by~The\-o\-rem~3.10 of \cite{MFsigma}. On the other hand, $J_2$ is generated by $H_2J^\sigma$ and $K_2J^\sigma$, which are subnormal subgroups of $G$, and hence Theorem 1.6.10 of \cite{lennox} yields that there is a subnormal subgroup $Y$ of $G$ such that $L_2M_2\leq Y\leq  H_2J^\sigma K_2J^\sigma$. Since $L_1M_1Y=YL_1M_1$, it follows from Theorem 3.15 of \cite{MFsigma} that $X=L_1M_1Y$ is $\sigma$-subnormal in $G$.

Finally, $$LM=L_1L_2M_1M_2=L_1M_1L_2M_2\leq X\leq L_1M_1H_2J^\sigma K_2J^\sigma=L_1H_2J^\sigma\cdot M_1K_2J^\sigma\subseteq HK$$ and the statement is proved.
\end{proof}

\begin{cor}
Let $\sigma=\{\sigma_i\,:\, i\in I\}$ be a partition of $\mathbb{P}$. Let $L$ and $M$ be subgroups of a locally finite group $G$. Let~$H$ and $K$ be $\sigma$-subnormal subgroups of $G$ which are minimal with respect to containing $L$ and $M$, respectively. If~\hbox{$LM=ML$,} then $HK=KH$ is $\sigma$-subnormal in~$G$.
\end{cor}
\begin{proof}
By Theorem \ref{altroteorem}, there is a $\sigma$-subnormal subgroup $X$ of $G$ such that $LM\leq X\leq HM$. Now, $X\cap H$ is $\sigma$-subnormal in $G$ and contains $M$, which means that \hbox{$X\cap H=H$.} Similarly, $X\cap K=K$, and so $X=HK=KH$. Finally, Theorem 3.15 of~\cite{MFsigma} shows that~$HK$ is $\sigma$-subnormal in $G$.
\end{proof}

\medskip

Finally, we deal with J.P. Williams' Join Theorem. Let $\sigma=\{\sigma_i\,:\, i\in I\}$ be a partition of $\mathbb{P}$ and let $A$ be a periodic abelian group. Recall that, for every set of primes $\pi$, $A_{\pi}$ denotes the $\pi$-component of $A$. Moreover, we let $D(A)$ be the largest divisible subgroup of $A$. We say that $A$ has {\it finite $\sigma$-join rank} if the following two conditions are satisfied:

\begin{itemize}

\item[(i)] there is no infinite subset $\tau$ of $\sigma$ such that $A_{\tau_i}/D(A_{\tau_i})$ has infinite rank for every $\tau_i\in\tau$;

\item[(ii)] there is no infinite subset $\tau$ of $\sigma$ such that $A_{\tau_i}/D(A_{\tau_i})$ has finite rank for every $\tau_i\in\tau$, but $A_{\rho}/D(A_{\rho})$ has infinite rank, where $\rho=\bigcup_{\tau_i\in\tau}\tau_i$.

\end{itemize}

\begin{theo}\label{williamstheo}
Let $\sigma=\{\sigma_i\,:\, i\in I\}$ be a partition of $\mathbb{P}$. Let $H$ and $K$ be locally finite groups.

\begin{itemize}

\item[\textnormal{(1)}]  If $H/H'\otimes K/K'$ has finite $\sigma$-join rank, then $\langle H,K\rangle$ is $\sigma$-subnormal in any locally finite group $G$ in which $H$ and $K$ can be $\sigma$-subnormally embedded.

\item[\textnormal{(2)}] If $H/H'\otimes K/K'$ has not finite $\sigma$-join rank, then there is a locally finite group $G$ containing $H$ and $K$ as subnormal subgroups, such that $\langle H,K\rangle$ is not $\sigma$-subnormal in $G$.

\end{itemize}
\end{theo}
\begin{proof}
(1)\quad Let $G$ be a locally finite group containing $H$ and $K$ as $\sigma$-subnormal subgroups. Let $J=\langle H,K\rangle$. Corollary \ref{corollaryperm2} yields that $J^\sigma$ is subnormal in $G$. Let $$J^\sigma =U_0\trianglelefteq U_1\trianglelefteq\ldots\trianglelefteq U_n=G$$ be the normal closure series of $J^\sigma$ in $G$. Since $J^\sigma$ is normalized by $J$, we have that every~$U_m$ is normalized by $J$ as well. Choose $0\leq\ell<n$. We claim that $J\, U_\ell/U_\ell$ is~\hbox{$\sigma$-sub}\-normal in $J\,U_{\ell+1}/U_\ell$. 

Since $H/H'\otimes K/K'$ has finite $\sigma$-join rank, we also have that $$HU_\ell/H'U_\ell\otimes KU_\ell/K'U_\ell\simeq H/H'(H\cap U_\ell)\otimes K/K'(K\cap U_\ell)$$ has finite $\sigma$-join rank. Thus, we can assume $U_\ell=\{1\}$; in particular, $J$ is $\sigma$-hypercentral. Now, it follows from~Lem\-ma~\ref{alternative} (and the remark before that statement) that there is a subset $\pi$ of~$\mathbb{P}$ such that $J=J_1\times J_2$, where $\pi(J_1)\subseteq\pi$, $\pi(J_2)\subseteq\mathbb{P}\setminus\pi$, $J_1$ is $\sigma$-nilpotent, $$H_2=H\cap J_2\trianglelefteq\trianglelefteq G\quad\textnormal{ and }\quad K_2=K\cap J_2\trianglelefteq\trianglelefteq G$$ (see also the proof of Theorem \ref{altroteorem}). Moreover, $J_1$ is generated by the $\sigma$-subnormal subgroups $H_1=H\cap J_1$  and $K_1=K\cap J_1$, and is therefore $\sigma$-subnormal in $G$ by~The\-o\-rem~3.6 of \cite{MFsigma}. Since $$H/H'\simeq(H_1/H_1')\times (H_2/H_2')\quad\textnormal{and}\quad K/K'\simeq(K_1/K_1')\times(K_2/K_2'),$$ we have that $$
\begin{array}{c}
\frac{H}{H'}\otimes\frac{K}{K'}\simeq \left(\frac{H_1}{H_1'}\otimes \frac{K_1}{K_1'}\right)\times\left(\frac{H_1}{H_1'}\otimes \frac{K_2}{K_2'}\right)\times \left(\frac{H_2}{H_2'}\otimes \frac{K_1}{K_1'}\right)\times\left(\frac{H_2}{H_2'}\otimes \frac{K_2}{K_2'}\right)
\end{array}
$$ and consequently that $H_2/H_2'\otimes K_2/K_2'$ has finite $\sigma$-join rank. Now, it easily follows from William's theorem (see \cite{lennox}, Theorem 5.1.1) that $J_2=\langle H_2,K_2\rangle$
is subnormal in~$JU_{\ell+1}$. Finally, Theorem 3.15 of \cite{MFsigma} shows that $J=J_1J_2$ is $\sigma$-subnormal in $JU_{\ell+1}$ and completes the proof of the claim.

In order to complete the proof of (1), we just need to observe that $$J=JU_0\trianglelefteq\trianglelefteq_\sigma JU_1\trianglelefteq\trianglelefteq_\sigma\ldots\trianglelefteq\trianglelefteq_\sigma JU_n=G$$ gives a finite $\sigma$-series connecting $J$ to $G$.

\medskip

(2) Let $\overline{H}=H/H'$ and $\overline{K}=K/K'$. It is clear that for each $\sigma_i\in\sigma$, the $\sigma_i$-component of  $\overline{H}\otimes\overline{K}$ is $\overline{H}_{\sigma_i}\otimes\overline{K}_{\sigma_i}$. This shows that there is an infinite set of primes $p_1,p_2,\ldots$ such that one of $\overline{H}$, $\overline{K}$ has a homomorphic image which is a direct sum of cyclic groups, $i$ of order $p_i$ for each $i$, and the other has a homomorphic image which is a direct sum of cyclic groups, one of order $p_i$ for each $i$. Moreover, we can (and we must) choose the primes $p_j$ in such a way that $p_m$ and $p_n$ belong to the same element of $\sigma$ if and only if $m=n$. Now, we apply the construction described at page 160 of \cite{lennox}, Case~(ii): here we must consider restricted direct sums instead of unrestricted ones, so the resulting group $G$ is locally finite. Our group can be described as $G=J\ltimes A$, where $A$ is abelian and $J=\langle H,K\rangle$. Moreover, $J$ contains a normal subgroup $N$ of $G$ and $J/N$ is $\sigma$-hypercentral. The argument at page 162 of \cite{lennox} shows that there is no infinite subset $\tau$ of $\sigma$ such that the largest $\tau$-nilpotent subgroup of $J/N$ is subnormal in $G/N$. Therefore $J/N$ cannot be $\sigma$-subnormal (see Lemma \ref{alternative}). The statement is proved.
\end{proof}

\begin{cor}
Let $\sigma=\{\sigma_i\,:\, i\in I\}$ be a partition of $\mathbb{P}$. Let $H_1,\ldots,H_n$ be $\sigma$-subnormal subgroups of a locally finite group $G$. If $H_i/H_i'\otimes H_j/H_j'$ has finite $\sigma$-join rank for all $i\neq j$, then $J=\langle H_1,\ldots,H_n\rangle$ is $\sigma$-subnormal in $G$.
\end{cor}
\begin{proof}
For each $1\leq i\leq n-1$, let $J_i=\langle H_1,\ldots, H_i\rangle$ and notice that $$J_i/J_i'\otimes H_{i+1}/H_{i+1}'$$ has finite $\sigma$-join rank. Thus, a repeated application of Theorem \ref{williamstheo} shows that $J$ is~\hbox{$\sigma$-sub}\-normal in $G$.
\end{proof}

 \medskip

 Recall that two subgroups $H$ and $K$ of a group $G$ are said to be {\it orthogonal} if the tensor product $H/H'\otimes K/K'$ is trivial.

\begin{cor}\label{roseblade}
Let $\sigma=\{\sigma_i\,:\, i\in I\}$ be a partition of $\mathbb{P}$. Let $H$ and $K$ be orthogonal~\hbox{$\sigma$-sub}\-normal subgroups of a locally finite group $G$. Then $\langle H,K\rangle$ is $\sigma$-subnormal in~$G$.
\end{cor}

\medskip

Although Roseblade's result shows that orthogonal subnormal subgroups permute, this is no longer true for orthogonal $\sigma$-subnormal subgroups, as shown by the examples just before Section \ref{sec3}.

\section{The class $\mathfrak{S}_\sigma$}\label{sec4}

In this section we deal with the properties of the class $\mathfrak{S}_\sigma$ of all locally finite groups in which the join of two $\sigma$-subnormal subgroups is $\sigma$-subnormal. It follows from~The\-o\-rem~\ref{corquestion2} that this class coincides with the class $\mathfrak{S}_{\operatorname{sn}\rightarrow\sigma}$ of all locally finite groups in which the join of two (actually, finitely many) subnormal subgroups is $\sigma$-subnormal. For each prime $p$, let $G_p$ be a metabelian $p$-group containing a subnormal subgroup of defect at least $p$, and let~$G$ be the direct product of all these groups. Since $G$ is metabelian, $G$ belongs to $\mathfrak{S}_\sigma$ by Theorem 3.32 of \cite{MFsigma}; but clearly $G$ does not belong to the class $\mathfrak{S}_{\operatorname{sn}\rightarrow\overline{\sigma}}^\infty$ of all locally finite groups in which the join of arbitrarily many subnormal subgroups is $\overline{\sigma}$-subnormal, where $\overline{\sigma}=\big\{\{2\},\{3\},\{5\},\ldots\big\}$. Needless to say, every locally finite primary group that does not belong to $\mathfrak{S}$ is an example of a group in $\mathfrak{S}_\sigma\setminus\mathfrak{S}$. The examples we presented at the beginning of this section and at the beginning of Section \ref{sec3} allow us to draw the following diagram of classes (there are no inclusions among the classes in the following diagram except those indicated in it, and these are all strict inclusions):

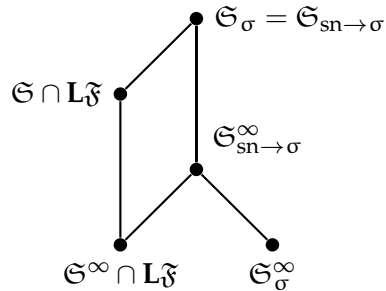
\begin{center}
\begin{tikzpicture}[thick, circ/.style={shape=circle, fill=black, inner sep=1pt, scale=1.5, draw, node contents=}]
\draw node (c1) at (1, 3) [circ, label=right:{$\mathfrak{S}_\sigma=\mathfrak{S}_{\operatorname{sn}\rightarrow\sigma}$}];
\draw node (c2) at (0, 2) [circ, label=left:{$\mathfrak{S}\cap{\bf L}\mathfrak{F}$}];
\draw node (c3) at (1, 1) [circ, label=above right:{$\mathfrak{S}_{\operatorname{sn}\rightarrow\sigma}^\infty$}];
\draw node (c4) at (0,0) [circ, label=below:{$\mathfrak{S}^\infty\cap{\bf L}\mathfrak{F}$}];
\draw node (c5) at (2, 0) [circ, label=below:{$\mathfrak{S}_\sigma^\infty$}];
\draw (c1) -- (c2) -- (c4) -- (c3) -- (c1) -- (c3) -- (c5);
\end{tikzpicture}
\end{center}

\noindent (here, ${\bf L}\mathfrak{F}$ denotes the class of all locally finite groups).

\medskip

The first goal of the section is to prove that $\mathfrak{M}\mathfrak{S}_\sigma\mathfrak{M}=\mathfrak{S}_\sigma$.  The fact that $\mathfrak{M}\mathfrak{S}_\sigma=\mathfrak{S}_\sigma$ can be established by an easy adaptation of the proof of~The\-o\-rem~\ref{sinftyminmax} using The\-o\-rem~\ref{lemmajoinfacile}.

\begin{theo}\label{sinftyminmaxfacile}
Let $\sigma=\{\sigma_i\,:\, i\in I\}$ be a partition of $\mathbb{P}$. Let $G$ be a locally finite group having a nor\-mal subgroup $N$ in $\mathfrak{M}$ such that $G/N\in\mathfrak{S}_\sigma$. Then $G$ is in $\mathfrak{S}_\sigma$.
\end{theo}

Now we turn to prove that $\mathfrak{S}_\sigma\mathfrak{M}=\mathfrak{S}_\sigma$. This is accomplished through some preliminary steps.


\begin{lem}\label{finiteextension}
Let $\sigma=\{\sigma_i\,:\, i\in I\}$ be a partition of $\mathbb{P}$. Let $G$ be a locally finite group having a normal subgroup $N\in\mathfrak{S}_\sigma$ of finite index. Then $G$ is in $\mathfrak{S}_\sigma$.
\end{lem}
\begin{proof}
By Theorem \ref{corquestion2}, it is enough to show that the join of a subnormal subgroup~$H$ of $G$ and a $\sigma$-subnormal subgroup $K$ of $G$ is always $\sigma$-subnormal in $G$. We use induction on the subnormal defect $m$ of~$H$ in $G$ (being the statement obvious if~\hbox{$m\leq1$)}. Let $S$ be a $\sigma$-embedded subgroup of $K$ that is subnormal in $G$ (see~Lem\-ma~\ref{alternative}). If $\langle H,S\rangle$ is~\hbox{$\sigma$-sub}\-normal in $G$, then $\langle H,K\rangle$ is~\hbox{$\sigma$-sub}\-normal in $G$ by~Lem\-ma~\ref{question2}. 

Now, let $A=S\cap N$, so $S/A$ is finite; let $n$ be the subnormal defect of $S$ in $G$. By~Co\-rol\-la\-ry 1.4.4 of \cite{lennox}, we can write $$H^S=L^A[H,_n\,S],$$ where $L$ is a subgroup generated by finitely many conjugates of $H$ under $S$. If we can prove that $H^S$ is $\sigma$-subnormal in $G$, then $\langle H,S\rangle=SH^S$ is $\sigma$-subnormal in $G$ by~The\-o\-rem~3.15 of \cite{MFsigma}, and we are done. But $[H,_nS]$ is a normal subgroup of $S$, so is subnormal in $G$; hence in order to prove that $H^S$ is $\sigma$-subnormal in $G$, we just have to prove that $L^A$ is $\sigma$-subnormal in $G$ (this is another application of~The\-o\-rem~3.15 of~\cite{MFsigma}).

By induction on $m$, we see that $L$ is $\sigma$-subnormal in $H^G$ and so in $G$. Let $T$ be a subnormal subgroup of $G$ that is $\sigma$-embedded in $L$ (see Lemma \ref{alternative}). If we can prove that $T^A$ is $\sigma$-subnormal in $G$, then $\langle A,T\rangle=AT^A$ is $\sigma$-subnormal in $G$ (by~The\-o\-rem~3.15 of~\cite{MFsigma}) and Lemma \ref{question2} shows that $\langle A,L\rangle$ is $\sigma$-subnormal in $G$; in particular,~$L^A$ is~\hbox{$\sigma$-sub}\-normal in $G$.

Let $\ell$ be the subnormal defect of $T$. Another application of Corollary 1.4.4 of \cite{lennox} gives $$A^T=B^{T\cap N}[A,_\ell\, T],$$ where $B$ is generated by finitely many conjugates of $A$. The sub\-group $$C=\langle B,\,T\cap N,\,[A,_\ell T]\rangle$$ is the join of finitely many subnormal subgroups of $N$, so $C$ is $\sigma$-subnormal in $N$ and hence in $G$. Since~$A^T$ is normal in $C$, it follows that $A^T$ is $\sigma$-subnormal in $G$. Finally,~\hbox{$\langle A,T\rangle=TA^T$} is $\sigma$-subnormal in $G$ by Theorem 3.15 of \cite{MFsigma} and consequently~$T^A$ is~\hbox{$\sigma$-sub}\-normal in $G$. The statement is proved.
\end{proof}

\begin{lem}\label{corextensionpezzi}
Let $\sigma=\{\sigma_i\,:\, i\in I\}$ be a partition of $\mathbb{P}$. Let $G$ be a locally finite group having normal subgroups $M\leq N$ such that $M$ and $G/N$ are $\sigma$-groups, while $N/M$ is in $\mathfrak{S}_\sigma$. Then~$G\in\mathfrak{S}_\sigma$.
\end{lem}
\begin{proof}
Let $H$ and $K$ be $\sigma$-subnormal subgroups of $G$, and let $J=\langle H,K\rangle$. The\-o\-rem~3.10 of \cite{MFsigma} shows that $J$ is~\hbox{$\sigma$-serial} in $G$ and Corollary 3.5 shows that $J$ is~\hbox{$\sigma$-sub}\-normal in $JM$. It is therefore possible to assume $M=\{1\}$. In this case,~\hbox{$N\cap H$} and~\hbox{$N\cap K$} are $\sigma$-embedded in $H$ and $K$, respectively. By hypothesis, the sub\-group~\hbox{$\langle N\cap H,\, N\cap K\rangle$} is $\sigma$-subnormal in $N$ and so in $G$. Thus, it follows from~Lem\-ma~\ref{question2} that $J$ is~\hbox{$\sigma$-sub}\-normal in $G$. The statement is proved.
\end{proof}

\medskip

We also need a way to reduce ourselves to deal with soluble groups. This is done by using an argument which is due to~Smith~\cite{Smith}.

\begin{theo}\label{solublereduction}
Let $\sigma=\{\sigma_i\,:\, i\in I\}$ be a partition of $\mathbb{P}$. Suppose that $\mathfrak{X}$ is a class of groups which is closed with respect to forming quotients and subnormal subgroups. If every soluble~\hbox{$\mathfrak{X}$-group} is in $\mathfrak{S}_\sigma$, then $\mathfrak{X}\leq\mathfrak{S}_\sigma$.
\end{theo}
\begin{proof}
Let $G$ be a group in $\mathfrak{X}$. It follows from Theorem \ref{corquestion2} that we only need to show that the join $J$ of two subnormal subgroups $H$ and $K$ is $\sigma$-subnormal. We argue by induction on the subnormal defect $m$ of $H$. Being the statement obvious for~\hbox{$m\leq1$,} we may assume $m\geq2$. Let $L$ be the first term of the normal closure series of $H$ in $G$ that is contained in $N_G(H)$. By the induction hypothesis, $J_0=\langle L,K\rangle$ is $\sigma$-subnormal in $G$. Let~\hbox{$P=P_L(K)$} be the permutizer of $K$ in $L$. Of course, $H^K=H^{PK}$ is a normal subgroup of $\langle H,PK\rangle=J_1$, which means that $J=H^KK$ is subnormal in $J_1$.

Now, by Theorem 1.6.8 of \cite{lennox} there is a positive integer $n$ with $L^{(n)}\leq P$, while by~The\-o\-rem~3.1.2 of \cite{lennox} there is a positive integer $\ell$ such that $$J_0^{(\ell)}\leq L^{(n)}K\subseteq PK\leq J_1\leq J_0.$$ Lemma \ref{alternative} shows that $J_0$ contains a $\sigma$-embedded subgroup $S$ that is subnormal in~$G$. It follows that $SJ^{(\ell)}_0/J^{(\ell)}_0$ is a soluble group in $\mathfrak{X}$ and hence $SJ^{(\ell)}_0/J^{(\ell)}_0$ is in $\mathfrak{S}_\sigma$. Repeated applications of Lemma \ref{corextensionpezzi} yield that $J_0/J_0^{(\ell)}$ is in $\mathfrak{S}_\sigma$.

Finally, $PK$ is subnormal in $J_0$ by Theorems 1.6.9 and 1.2.5 of \cite{lennox}; consequently, $J_1=\langle H,PK\rangle$ is $\sigma$-subnormal in $J_0$, and $$J\trianglelefteq\trianglelefteq J_1\trianglelefteq\trianglelefteq_\sigma J_0\trianglelefteq\trianglelefteq_\sigma G.$$ The statement is proved.
\end{proof}

\medskip

We can now prove the desired analog of Theorem \ref{sinftyminmaxfacile}.

\begin{theo}\label{theosopra}
Let $\sigma=\{\sigma_i\,:\, i\in I\}$ be a partition of $\mathbb{P}$. Let $G$ be a locally finite group having a normal subgroup $N\in\mathfrak{S}_\sigma$ such that $G/N$ is in $\mathfrak{M}$. Then $G$ is in $\mathfrak{S}_\sigma$.
\end{theo}
\begin{proof}
By Theorem \ref{solublereduction}, we may assume $G$ is soluble. Moreover, it is certainly possible to assume $G/N$ is in $\operatorname{Min}\textnormal{-}\operatorname{sn}$ or in $\operatorname{Max}\textnormal{-}\operatorname{sn}$. If $G/N$ satisfies the maximal condition on subnormal subgroups, then $G/N$ is finite and the results follows from~Lem\-ma~\ref{finiteextension}. Assume $G/N$ satisfies the minimal condition on subnormal subgroups. In this case, $G/N$ is \v Cernikov, and a combination of   Lemma \ref{finiteextension} and Lemma \ref{corextensionpezzi} completes the proof of the statement.
\end{proof}

\medskip

A careful look at the above proofs shows that in Theorems \ref{sinftyminmaxfacile} and \ref{theosopra} we may replace $\mathfrak{M}$ by any class of groups $\mathfrak{X}\leq\mathfrak{S}_\sigma$ which is closed with respect to forming subnormal subgroups, quotients and is such that every hyper-$\sigma$ $\mathfrak{X}$-group is $\sigma$-soluble (actually, we just need that every $\sigma$-hypercentral $\mathfrak{X}$-group is $\sigma$-nilpotent).

\medskip

We now wish to replace the class $\mathfrak{M}$ by the class of all locally finite groups with finite rank. First, we extrapolate the argument used in the proof of Theorem 3.5.1 of~\cite{lennox} to prove the following lemma.

\begin{lem}\label{lemmarango}
Let $G=JA$ be a group, where $A$ is normal abelian of finite rank $r$ and $J$ is generated by two subnormal subgroups $H$ and $K$ of $G$. Then $J$ is subnormal in $G$.
\end{lem}
\begin{proof}
Since $J\cap A$ is normal in $G$, we may assume $G=J\ltimes A$. Similarly, we may assume $C_J(A)=\{1\}$. Now, a theorem of Philip Hall yields that $H$ and $K$ are nilpotent; actually, if $m$ and $n$ are the defects of $H$ and $K$, respectively, then $H$ and $K$ are nilpotent of classes at most $m-1$ and $n-1$, respectively; in particular, $G$ is locally nilpotent. 

Let $T$ be the periodic part of $A$. It follows from Lemma 6.37 of \cite{Rob72} that $A/T\leq\zeta_r(G/T)$ for some $r$, so $JT$ is subnormal in $G$, and we may assume $A=T$ is periodic.  

We prove that $[A,_\ell J]=\{1\}$ for some $\ell$ depending only on $n$, $m$ and $r$. To this aim, we need some further reductions. Since $A$ is periodic, it is possible to assume $A$ is a~\hbox{$p$-group} for some prime $p$. Let $H_0$ and $K_0$ be finitely generated subgroups of $H$ and~$K$, respectively, and let $J_0=\langle H_0,K_0\rangle$. Since $H$ and $K$ are nilpotent, their respective subgroups $H_0$ and~$K_0$ are subnormal in $G$ of defect bounded by a function of $m$ and $n$ only. Let~\hbox{$a\in A$} and $B=\langle a\rangle^{J_0}$, so $B$ is a finite $p$-group. It is thus possible to replace~$G$ by $\big(J_0\ltimes B\big)/C_{J_0}(B)$; in particular,~$G$ is a finite~\hbox{$p$-group} (because such is $J_0/C_{J_0}(B)$) of rank bounded by a function of $r$ (see~\cite{lennox},~Lem\-ma~6.1.10). 

Now, the nilpotency class of $J$ is bounded by a function of $m$, $n$ and $r$ only (see~The\-o\-rem~3.3.1 of~\cite{lennox}) and consequently Theorem 3.3.2 of \cite{lennox} shows that $$[A,_\ell\, J]\leq [A,_mH][A,_nK]=\{1\}$$ for some $\ell$ depending only on $n$, $m$ and $r$. The statement is proved.
\end{proof}

\begin{theo}\label{theorankunico}
Let $\sigma=\{\sigma_i\,:\, i\in I\}$ be a partition of $\mathbb{P}$. Let $G$ be a locally finite group having a normal subgroup $N$ of finite rank such that $G/N\in\mathfrak{S}_\sigma$. Then $G$ is in $\mathfrak{S}_\sigma$.
\end{theo}
\begin{proof}
It follows from Theorem \ref{solublereduction} that we may assume $G$ is soluble. We work by induction on the derived length of $N$. By Theorem \ref{corquestion2}, we only need to show that the join $J$ of two subnormal subgroups $H$ and $K$ of $G$ is $\sigma$-subnormal in $G$. 

Let $A$ be the last non-trivial term of the derived series of $N$. By induction hypothesis,~$JA$ is $\sigma$-subnormal in~$G$. On the other hand,~Lem\-ma~\ref{lemmarango} shows that $J$ is subnormal in $G$ and completes the proof of the statement.
\end{proof}

\begin{cor}
Let $\sigma=\{\sigma_i\,:\, i\in I\}$ be a partition of $\mathbb{P}$. If $G$ is a locally finite group of finite rank, then $G\in\mathfrak{S}_\sigma$.
\end{cor}

\medskip

The situation is much less clear for a locally finite group $G$ that is the extension of a normal subgroup~$N\in\mathfrak{S}_\sigma$ (for a certain partition $\sigma$ of $\mathbb{P}$) by a group of finite rank. Of course, by Theorem~\ref{solublereduction} we may assume $G$ is soluble. Then working by induction on the derived length of $G$, we may easily reduce to the case $G=J\ltimes A$, where $A$ is abelian and $J$ is generated by two nilpotent subnormal subgroups $H$ and $K$ of $G$ (see~The\-o\-rem~\ref{corquestion2}). It follows that $G$ is locally nilpotent (actually, $G$ is a Baer group), so $G$ is the direct product of its primary components. The goal is to prove that $J$ is~\hbox{$\sigma$-sub}\-normal in $G$, but unfortunately this is as far as we can go.

\begin{quest}\label{questionpos}
Let $\sigma=\{\sigma_i\,:\, i\in I\}$ be a partition of $\mathbb{P}$. Is it true that any locally finite extension of a group in $\mathfrak{S}_\sigma$ by a group of finite rank is in $\mathfrak{S}_\sigma$?
\end{quest}

The above problem is perhaps connected with the lack of information concerning an extension of a group in $\mathfrak{S}$ by a group of finite rank; and we actually conjecture the answer to Question \ref{questionpos} is negative. On the other hand, if we strengthen the hypotheses on the normal subgroup $N$, we can give a positive answer (see also~\cite{Smith},~The\-o\-rem~2.7). In what follows, denote by $\mathfrak{Y}$ the class of all groups $G$ having a normal subgroup $N$ such that $G/N$ and $N'$ have finite rank.

\begin{theo}\label{dacitare}
Let $\sigma=\{\sigma_i\,:\, i\in I\}$ be a partition of $\mathbb{P}$. Let $G$ be a locally finite group having a normal subgroup $N\in\mathfrak{S}_{\operatorname{sn}\rightarrow\sigma}^\infty$ such that $G/N\in\mathfrak{Y}$. Then $G$ is in $\mathfrak{S}_{\sigma}$.
\end{theo}
\begin{proof}
By Theorem \ref{solublereduction}, we may assume $G$ is soluble. Let $J$ be the subgroup generated by two subnormal subgroups $H$ and $K$ of $G$. The goal is to prove that $J$ is $\sigma$-subnormal in $G$ (see Theorem \ref{corquestion2}). The\-o\-rem~\ref{lemmajoinfacile} yields that~\hbox{$J\cap N$} is~\hbox{$\sigma$-sub}\-normal in $N$ and so in~$G$. Then there is a subnormal subgroup~$S$ of~$G$ that is $\sigma$-embedded in~\hbox{$J\cap N$} and~\hbox{$J$-in}\-va\-riant. Let $$S=S_0\trianglelefteq S_1\trianglelefteq\ldots\trianglelefteq S_n=N$$ be the normal closure series of $S$ in $N$. Choose $0\leq i<n$. We have that~$JS_{i}/S_i$ is~\hbox{$\sigma$-sub}\-normal in $JL/S_i$, where $L/S_i$ is the normal closure of $(J\cap N)S_{i}/S_i$ in~$JS_{i+1}/S_i$ (see~\cite{MFsigma},~Lem\-ma~2.6, Corollary 3.5 and Theorem 3.10). Thus, $JL/L\in\mathfrak{Y}$ and $JL\cap N_1=L$, where $N_1=JS_{i+1}\cap N=(J\cap N)S_{i+1}$. Thus, $$JS_{i+1}/L=JL/L\ltimes N_{1}/L.$$ On the other hand, $JL/L$ is generated by two subnormal subgroups $H/L$ and $K/L$ of~$JS_{i+1}/L$, and $N_1/L$ is soluble, so using the already mentioned theorem of Philip~Hall we can find a subgroup $M/L$ of $JL/L$ which is normal in $JS_{i+1}/L$ and such that $HM/M$ and $KM/M$ are nilpotent. It follows from Corollary 4.8 of \cite{Smith} that $JL/M$ is nilpotent and hence $JL/M$ is subnormal in $JS_{i+1}/M$ by Lemma 4.10 of \cite{Smith}. Therefore $JS_i$ is~\hbox{$\sigma$-sub}\-nor\-mal in~$JS_{i+1}$. Combining all these $\sigma$-series, we find that $J=JS_0$ is~\hbox{$\sigma$-sub}\-normal in $G$. The statement is proved.
\end{proof}

\begin{cor}\label{corabelianinf}
Let $\sigma=\{\sigma_i\,:\, i\in I\}$ be a partition of $\mathbb{P}$. Let $G$ be a locally finite group having a normal subgroup $N\in\mathfrak{S}_{\operatorname{sn}\rightarrow\sigma}^\infty$ such that $G/N$ is abelian. Then $G$ is in $\mathfrak{S}_\sigma$.
\end{cor}

\medskip

Corollary \ref{corabelianinf} can be extended in another direction. In what follows, let $\mathfrak{T}$ be the class of groups in which normality is a transitive relation.

\begin{theo}\label{corabelianinftheo}
Let $\sigma=\{\sigma_i\,:\, i\in I\}$ be a partition of $\mathbb{P}$. Let $G$ be a locally finite group having a normal subgroup $N\in\mathfrak{S}_{\operatorname{sn}\rightarrow\sigma}^\infty$ such that $G/N\in\mathfrak{T}$. Then $G$ is in $\mathfrak{S}_\sigma$.
\end{theo}
\begin{proof}
Let $J$ be the subgroup generated by two subnormal subgroups $H$ and $K$ of $G$. The goal is to prove that $J$ is $\sigma$-subnormal in $G$ (see Theorem \ref{corquestion2}). The\-o\-rem~\ref{lemmajoinfacile} yields that~\hbox{$J\cap N$} is~\hbox{$\sigma$-sub}\-normal in $N$ and so in~$G$. Then there is a subnormal subgroup~$S$ of~$G$ that is $\sigma$-embedded in~\hbox{$J\cap N$} and~\hbox{$J$-in}\-va\-riant. Let $$S=S_0\trianglelefteq S_1\trianglelefteq\ldots\trianglelefteq S_n=N$$ be the normal closure series of $S$ in $N$. Choose $0\leq i<n$. We have that~$JS_{i}/S_i$ is~\hbox{$\sigma$-sub}\-normal in $JL/S_i$, where $L/S_i$ is the normal closure of $(J\cap N)S_{i}/S_i$ in~$JS_{i+1}/S_i$ (see~\cite{MFsigma},~Lem\-ma~2.6, Corollary 3.5 and Theorem 3.10). Since $JN/N$ is normal in~$G/N$, it follows that $JL/L$ is in $\mathfrak{T}$. This means that $JL/L=HL/L\cdot KL/L$, so $JL/L$ is $\sigma$-subnormal in $JS_{i+1}/S_i$ by Theorem 3.15 of \cite{MFsigma}. Therefore $JS_i$ is~\hbox{$\sigma$-sub}\-nor\-mal in~$JS_{i+1}$. Combining all these $\sigma$-series, we find that $J=JS_0$ is~\hbox{$\sigma$-sub}\-normal in $G$. The statement is proved.~\end{proof}

\medskip

Theorem 3.32 of \cite{MFsigma} (or Corollary \ref{corabelianinf} of this paper) shows that the class of locally finite, metabelian groups is contained in $\mathfrak{S}_\sigma$. Thus, the example constructed in Theorem 5.1 of \cite{Smith} (which is soluble of length $3$) shows that we cannot replace $\mathfrak{S}_{\operatorname{sn}\rightarrow\sigma}^\infty$ by $\mathfrak{S}_\sigma$ in the three above results (we note in passing that the same example is easily seen to be the product of two normal metabelian groups). Moreover the example just before Theorem 3.10 of \cite{MFsigma} shows that in these results we cannot replace ``$G$ is in $\mathfrak{S}_\sigma$'' by ``$G$ is in $\mathfrak{S}_\sigma^\infty$''. 

We also remark that a combination of Theorem 3.4.3 of \cite{lennox} and our Lemma \ref{question2} yields that the join $J$ of two $\sigma$-subnormal subgroups $H$ and $K$ of a locally finite group is $\sigma$-subnormal provided that $J'$ has finite rank. Here we cannot replace ``finite rank'' by the class $\mathfrak{Y}$ defined above because of the example in \cite{Smith} ($J$ is metabelian in this example).

Furthermore, the quoted example from \cite{Smith} shows that we cannot swap places of abelianity and $\mathfrak{S}_{\operatorname{sn}\rightarrow\sigma}^\infty$ in the statement of Corollary \ref{corabelianinf}.

\medskip

Finally, we briefly deal with $\sigma$-subnormal subgroups in case the partition $\sigma$ is finite. Here the main result is the following one.

\begin{theo}\label{finalone}
Let $\sigma=\{\sigma_i\,:\, i\in I\}$ be a finite partition of $\mathbb{P}$. If $G$ is a locally finite group, then $G$ belongs to $\mathfrak{S}_\sigma$.
\end{theo}
\begin{proof}
This is a consequence of Corollary 3.31 of \cite{MFsigma} and our Theorem \ref{solublereduction}.
\end{proof}

\section{Declarations}

\noindent{\bf Conflict of interest}\quad The authors declare that they have no conflict of interest.

\smallskip

\noindent{\bf Data availability}\quad Data sharing is not applicable to this article as no new data were created or analyzed in this study.

\bigskip\bigskip\bigskip

\renewcommand{\bibsection}{\begin{flushright}\Large
{
REFERENCES}\\
\rule{8cm}{0.4pt}\\[0.8cm]
\end{flushright}}

\bigskip\bigskip


\begin{flushleft}
\rule{8cm}{0.4pt}\\
\end{flushleft}

{
\sloppy
\noindent
Maria Ferrara

\noindent
Dipartimento di Matematica e Fisica

\noindent
Università degli Studi della Campania  ``Luigi Vanvitelli''

\noindent
viale Lincoln 5, Caserta (Italy)

\noindent
e-mail: maria.ferrara1@unicampania.it
}

\bigskip
\bigskip

{
\sloppy
\noindent
Marco Trombetti

\noindent 
Dipartimento di Matematica e Applicazioni ``Renato Caccioppoli''

\noindent
Università degli Studi di Napoli Federico II

\noindent
Complesso Universitario Monte S. Angelo

\noindent
Via Cintia, Napoli (Italy)

\noindent
e-mail: marco.trombetti@unina.it 

}


\end{document}